\documentclass[12pt]{article}
\pdfoutput=1
\usepackage{epsfig,amsmath,amsfonts,amssymb,latexsym,color,amsthm,hhline,multirow,subfigure,graphicx}

\newtheorem{theorem}{Theorem}[section]

\newcommand{\E}{{\mathbf E}}

\newcommand {\PP}{{\mathbb P}}
\newcommand{\sss}{\scriptscriptstyle}

\begin{document}

\title{First passage percolation on $\mathbb{Z}^2$ -- a simulation study}\parskip=5pt plus1pt minus1pt \parindent=0pt
\author{Sven Erick Alm\thanks{Department of Mathematics, Uppsala University, 751 06 Uppsala, Sweden; {\tt sea@math.uu.se}}\and Maria Deijfen\thanks{Department of Mathematics, Stockholm University, 106 91 Stockholm, Sweden; {\tt mia@math.su.se}}}
\date{December 2014}
\maketitle

\begin{abstract}
\noindent First passage percolation on $\mathbb{Z}^2$ is a model for describing the spread of an infection on the sites of the square lattice. The infection is spread via nearest neighbor sites and the time dynamic is specified by random passage times attached to the edges. In this paper, the speed of the growth and the shape of the infected set is studied by aid of large-scale computer simulations, with focus on continuous passage time distributions. It is found that the most important quantity for determining the value of the time constant, which indicates the inverse asymptotic speed of the growth, is $\E[\min\{\tau_1,\ldots,\tau_4\}]$, where $\tau_1,\ldots,\tau_4$ are i.i.d.\ passage time variables. The relation is linear for a large class of passage time distributions. Furthermore, the directional time constants are seen to be increasing when moving from the axis towards the diagonal, so that the limiting shape is contained in a circle with radius defined by the speed along the axes. The shape comes closer to the circle for distributions with larger variability.

\vspace{0.5cm}


\end{abstract}

\section{Introduction}

First passage percolation was introduced by Hammersley and Welsh in \cite{HW65} and can be viewed as a model for the spread of an infecion on a graph structure. The term 'infection' can have different meanings depending on the context, and may refer e.g.\ to a fluid spreading through a porous medium, information transfer or a virus spreading in a structured population. The original and most studied case is when the underlying graph is the $\mathbb{Z}^d$ lattice. Here we restrict to $d=2$, and investigate the speed and shape of the growing infected set by aid of computer simulations. The theoretical results involve quantities that have turned out to be very difficult to characterize analytically and simulation studies may therefore be of interest to shed light on qualitative behavior and as guidelines for further theoretical investigations. The total CPU time for our simulations amounts to more than 19 years, which make them the most extensive ones so far for the model.

To specify the model, let $\mathbb{Z}^2$ denote the set of integer points in the plane and write $\mathbb{E}$ for the set of neareast neighbor edges of $\mathbb{Z}^2$. Each edge $e\in\mathbb{E}$ is equipped with a non-negative random variable $\tau(e)$ interpreted as the time it takes for the infection to traverse the edge $e$. A path is a sequence of connected edges $e_1,\ldots,e_n\subset \mathbb{E}$ and the passage time for a path $\Gamma$ is given by
$$
T(\Gamma)=\sum_{e\in\Gamma}\tau(e).
$$
The passage time from a set $A\subset \mathbb{Z}^2$ to a set $B\subset \mathbb{Z}^2$ is defined as
$$
T(A,B)=\inf\{T(\Gamma):\Gamma\mbox{ is a path starting in $A$ and ending in $B$}\}.
$$
When $A$ and $B$ consist of single sites, with $A=\{x\}$ and $B=\{y\}$, we write $T(A,B)=T(x,y)$.

The interest revolves around dynamics and asymptotics for the growth of the infection, starting at time 0 with the origin infected and all other sites of $\mathbb{Z}^2$ uninfected. In order to be able to analyze this, one has to make some assumptions on the joint distribution of the set of edge passage times $\{\tau(e)\}_{e\in\mathbb{E}}$. Usually, the passage times are assumed to be i.i.d., but some results requiring only stationarity also appear in the literature, see e.g.\ \cite{B90,HM95}. Here we shall throughout assume that the passage times are i.i.d.\ with common distribution function $F$.

Existing results basically fall in three categories: (i) results on the growth in a fixed direction, (ii) results concerning the growth of the whole set of infected sites, and (iii) results on properties of infinite time-minimizing paths (geodesics). Our purpose here is to explore results of type (i) and (ii) by aid of simulations. We will not be concerned with results in category (iii), but refer e.g.\ to \cite{LN96,N95} for early work and conjectures, and to \cite{DamHan12,Hoff08} for more recent results and proofs of some of the conjectures. We also remark that such results are related to properties of growth processes involving several competing infection types; see \cite{DH08} and the references therein.

To characterize the growth in a fixed direction, first note that the passage times defined above are subadditive in the sense that
$$
T(x,y)\leq T(x,z)+T(z,y)\quad \mbox{for all }x,y,z\in\mathbb{Z}^2.
$$
Writing $\mathbf{n}=(n,0)$, this means that, under sufficient moment conditions on the distribution $F$ of the edge passage times, the family $\{T(\mathbf{m},\mathbf{n})\}_{m,n\in\mathbb{N}}$ fulfills the conditions of Liggett's version \cite{L85} of Kingman's subadditive ergodic theorem \cite{K68}, which asserts roughly that $T(\mathbf{0},\mathbf{n})$ grows linearly in $n$ as $n\to\infty$. In order for the theorem to be applicable, we need that $\E[T(\mathbf{0},\mathbf{n})]<\infty$ for all $n\geq 1$ and, by \cite[Theorem 4.6]{SW78}, this is the case if and only if
\begin{equation}\label{eq:minfin}
\E[\min\{\tau_1,\ldots, \tau_4\}]<\infty,
\end{equation}
where $\tau_1,\ldots,\tau_4$ are i.i.d.\ with distribution $F$. Hence (\ref{eq:minfin}) is the required moment condition, and the conclusion of Liggett's theorem then is that
$$
\frac{T(\mathbf{0},\mathbf{n})}{n}\to \mu \quad\mbox{ a.s. and in }L^1,
$$
where
\begin{equation}\label{eq:mu_def}
\mu=\inf_{n\geq 1}\frac{\E[T(\mathbf{0},\mathbf{n})]}{n}<\infty.
\end{equation}
The constant $\mu$ is known as the time constant. Its inverse value gives the asymptotic speed of the growth along the axes and the fact that $\mu<\infty$ implies that the infection grows at least linearly in time. In \cite{Kes86}, it is shown that $\mu>0$ if and only if
\begin{equation}\label{eq:mu_pos}
F(0)<1/2.
\end{equation}
Hence the infection grows linearly in time when $F(0)<1/2$ and faster than linearly when $F(0)\geq 1/2$. We remark that in dimensions $d\geq 3$, the bound 1/2 is replaced by the critical value of standard bond percolation on $\mathbb{Z}^d$.

A fundamental problem in first passage percolation is to determine the time constant $\mu$ and its relation to the edge passage time distribution $F$. This has turned out to be a difficult task: the value of $\mu$ is so far not known for any non-degenerate distribution $F$. Existing results and bounds are quite weak and are described in Section 2. One of the aims of this paper is to investigate by simulations how $\mu$ is related to various characteristics of $F$.

In addition to studying asymptotics for passage times along the axes, it is of course also of interest to study passage times in other directions. For $x\in\mathbb{R}^2$, we interpret $T(\mathbf{0},x)$ as the time when the closest point of $x$ in $\mathbb{Z}^2$ is infected. With $\hat{x}=x/|x|$, the subadditivity arguments outlined above can then be adapted to establish the existence of a directional time constant $\mu(\hat{x})$ such that $\lim T(\mathbf{0},n\hat{x})/n=\mu(\hat{x})$ a.s.\ (some extra work is required for $x$ with non-rational coordinates). Here, $\mu(\mathbf{1})=\mu$, where $\mu$ is the time constant along the axes. An interesting question is of course if and how $\mu(\hat{x})$ changes as the direction of $\hat{x}$ varies. This is closely related to the nature of the asymptotic shape of the infected set, which is our next topic.

We now proceed to study the growth of the whole infected set. To this end, write $B(t)=\{x\in\mathbb{Z}^2:T(\mathbf{0},x)\leq t\}$, let $C_x$ denote a unit cube centered at $x\in\mathbb{R}^2$ and define $\bar{B}(t)=\cup_{x\in B(t)} C_x$. The main result is a shape theorem, which asserts that $\bar{B}(t)/t$ converges almost surely to a deterministic shape $A$. The following version is due to Cox and Durrett \cite{CD81}. Kesten \cite{Kes86} has proved that the conclusion is true also in dimension $d\geq 3$ under stronger moment conditions.

\begin{theorem}\label{th:shape} Let $\tau_1,\dots, \tau_4$ be i.i.d.\ with distribution $F$ and assume that
\begin{equation}\label{eq:shape_cond}
\E[\min\{\tau_1^2,\ldots,\tau_4^2\}]<\infty.
\end{equation}
\begin{itemize}
\item[\rm{(a)}] If $F(0)<1/2$, then there exists a non-random compact convex set $A\subset \mathbb{R}^2$ with non-empty interior such that, for all $\varepsilon>0$, we have a.s.\ that
$$
(1-\varepsilon)A\subset \frac{\bar{B}(t)}{t}\subset (1+\varepsilon) A\quad\textrm{for large }t.
$$
\item[\rm{(b)}] If $F(0)\geq 1/2$, then for all $m\geq 0$ we have a.s.\ that
$$
\frac{\bar{B}(t)}{t}\supset \{x\in\mathbb{R}^2: |x|\leq m\}\quad\textrm{for large }t.
$$
\end{itemize}
If (\ref{eq:shape_cond}) fails, then $\limsup_{x:|x|\to\infty}\frac{T(\mathbf{0},x)}{|x|}\to\infty$ a.s.
\end{theorem}

The last statement implies that (\ref{eq:shape_cond}) is a necessary condition for the linear growth of $B(t)$, since otherwise the infected region will have uninfected points in its interior that are visible on a linear scale. We refer to Section \ref{sim_shape} for pictures illustrating this. Furthermore, part (b) of the theorem asserts that, if $F(0)\geq 1/2$, then the scaled infected set contains any bounded region of $\mathbb{R}^2$ eventually a.s.\ We will be interested in the case when $F(0)<1/2$ and (\ref{eq:shape_cond}) holds, so that the growth is linear and the scaled infected set converges to a compact shape $A$. Apart from the symmetries inherited from the $\mathbb{Z}^2$ lattice, not much is known about $A$, and the difficulties with characterizing $A$ analytically basically stem from lattice effects caused by $\mathbb{Z}^2$. We will instead study $A$ by aid of simulations. Questions that we will try to answer include: \emph{How does the shape deviate from a circle? How is it affected by the passage time distribution?  What about the fluctuations of the shape around its mean?}

The rest of the paper is organized so that Section 2 contains a survey of known results on the time constant and the asymptotic shape, in Section 3 the results of the simulations are described and Section 4 contains conclusions and directions for further work.

\section{Theoretical results}

We begin with a short survey of known results that are relevant to our studies. We also define the quantities that will be investigated and describe the aims of the simulations in more detail. The material is organized in two sections, one on growth in a fixed direction and one on the asymptotic shape.

\subsection{Growth in a fixed direction}

As mentioned in the previous section, the time constant $\mu$ has turned out to be difficult to determine. As for rigorous bounds, it follows from the characterization (\ref{eq:mu_def}) that $\mu\leq \E[\tau(e)]$, and in \cite{HW65} it is shown that this inequality is strict for all non-degenerate distributions $F$. Methods for computing better rigorous upper bounds can be found in \cite{AP02,SW78}, and lower bounds are derived in \cite{AP02,J81}. Both the upper and the lower bounds in general require computer assistance for obtaining numerical values. The precision varies with the distribution $F$, but the bounds tend to be quite crude. For instance, the methods in \cite{AP02} give $0.30\leq \mu \leq 0.503$ for Exp(1) passage times and $0.243\leq\mu\leq 0.403$ for $U(0,1)$ passage times, while the estimates (from computer simulations) given in \cite{AP02} are 0.402 and 0.312, respectively. 

Write $\mu^{\sss F}$ for the time constant associated with the edge passage time distribution $F$. An interesting question is how $\mu^{\sss F}$ behaves as a functional of $F$. Cox and Kesten \cite{CK81} show that $\mu^{\sss F}$ is continuous in $F$ in the sense that, if $F_n$ converges weakly to $F$, then $\mu^{\sss F_n}\to \mu^{\sss F}$. Another natural property, confirmed in \cite{SW78}, is that, if $F(x)\leq \widetilde{F}(x)$ for all $x\in\mathbb{R}$, then
\begin{equation}\label{eq:mu_ineq}
\mu^{\sss \widetilde{F}}\leq \mu^{\sss F}.
\end{equation}
In \cite{vdBK93}, it is shown that (\ref{eq:mu_ineq}) holds also under the weaker assumption that
$$
\int\phi(x)d\widetilde{F}(x)\leq \int\phi(x)dF(x)
$$
for every concave increasing function $\phi$ for which the above integrals are absolutely convergent. If this holds, then $\widetilde{F}$ is said to be more variable than $F$. Furthermore, if $\widetilde{F}$ is more variable than $F$ and $\widetilde{F}\neq F$, then $\widetilde{F}$ is said to be strictly more variable than $F$. The main result in \cite{vdBK93} is that the inequality in (\ref{eq:mu_ineq}) is strict if $\widetilde{F}$ is strictly more variable than $F$ and if (\ref{eq:mu_pos}) holds (an additional assumption for distributions with inf supp$(F)>0$ was later removed in \cite{M02}, where supp($F$) denotes the support of $F$).


One way of transforming a given distribution into a strictly more variable one is by pushing mass away from some point $\xi$; for further details see the criterion by Karlin and Novikoff in \cite{SD83}. A simple example is when $F$ is uniform on $[a,b]$ and $\widetilde{F}$ is uniform on $[a-\varepsilon_1,b+\varepsilon_2]$ for $\varepsilon_1,\varepsilon_2>0$. This example demonstrates that two distributions with the same mean do not necessarily give rise to the same time constant. For the particular case of uniform distributions centered at a given mean, the time constant is instead strictly decreasing as the variance increases. An interesting question is of course what properties of $F$ that are important in determining $\mu^{\sss F}$ in general. This largely remains unknown, and is one of the questions that we will investigate here. \emph{Is it true in general that, for distributions with the same mean, the time constant is decreasing in the variance? Can we find another quantity than the mean that determines the time constant?}

We continue by exploring passage times to lines rather than to single points. To this end, write $L_n$ for the vertical line that crosses the $x$-axis at $\mathbf{n}$. It turns out that the passage time $T(\mathbf{0},L_n)$ to the line $L_n$ behaves asymptotically the same as $T(\mathbf{0},\mathbf{n})$, that is,
$$
\frac{T(\mathbf{0},L_n)}{n}\to \mu \quad\mbox{ a.s. and in }L^1.
$$
A first proof of this appeared in \cite{WiRe78}, but it can also be derived from the asymptotic shape result, stated below as Theorem \ref{th:shape}; see \cite[pp.166-167]{Kes86}. The observation that $\mu$ appears also as the limit of the scaled passage times to the line $L_n$ has advantages when estimating $\mu$ by aid of simulations: The passage time of the infection to the single point $\mathbf{n}$ may be substantially increased if all edges incident to $\mathbf{n}$ happen to have very large passage times, and keeping track of the whole growth process during this time requires a lot of computer power. This problem arises in particular for heavy-tailed edge passage time distributions $F$, but is avoided by working with passage times to lines rather than points.

\begin{figure}
\begin{center}
\mbox{\includegraphics[height=4cm]{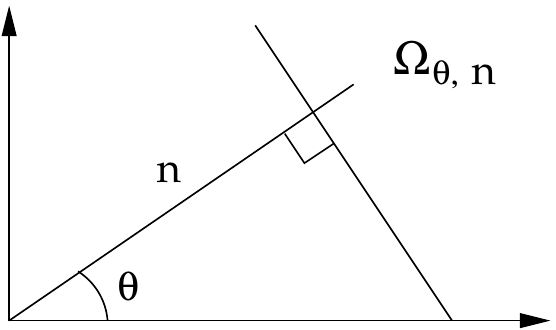}}
\end{center}
\caption{The region $\Omega_{\theta,n}$.}
\label{fig:omega}
\end{figure}

To define passage times to lines in an arbitrary direction, fix an angle $\theta\in[0,\pi/4]$ and write $\alpha=\tan \theta$. Furthermore, let $\Omega_{\theta,n}=\{(x_1,x_2)\in\mathbb{R}^2:x_1+\alpha x_2\geq n/\cos\theta\}$, that is, $\Omega_{\theta,n}$ is the part of $\mathbb{R}^2$ to the right of the line that is orthogonal to the line $y=\alpha x$ and intersects this line at distance $n$ from the origin; see Figure \ref{fig:omega}. Assuming (\ref{eq:minfin}), it follows from subadditivity that
$$
\frac{T(\mathbf{0},\Omega_{\theta,n})}{n}\to \mu_\theta \quad\textrm{a.s. and in }L^1,
$$
where
$$
\mu_\theta:=\inf_{n\geq 1}\frac{\E[T(\mathbf{0},\Omega_{\theta,n})]}{n}.
$$
Here $\mu_\theta$ is the time constant in the direction $\theta$ based on passage times to lines. The arguments for $\theta=0$ (referred to above) can easily be adapted to show that the constant coincides with the one based on passage times to points, that is, $\mu_\theta=\mu(\hat{x}_\theta)$, where $x_\theta=(1,\alpha)$.

Clearly $\mu_0=\mu$, where $\mu$ is the time constant along the axes, described above. An interesting question is of course if and how $\mu_\theta$ changes as $\theta$ increases. However, just as for $\theta=0$, obtaining information about $\mu_\theta$ for $\theta>0$ is difficult. The methods in \cite{AP02,J81} for calculating rigorous bounds are applicable in any direction, but there is no systematic study on how the resulting bounds depend on the direction: As mentioned above, the bounds tend to be quite crude and this means that they will most likely not allow for rigorous conclusions on how $\mu_\theta$ changes with $\theta$. We will instead study this by aid of simulations. The behavior of $\mu_\theta$ as a function of $\theta$ is closely related to the nature of the asymptotic shape of the infected set, which is our next topic.

\subsection{Asymptotic shape}\label{sec:AS}

Recall that $\mu(x)=\lim_nT(\mathbf{0},nx)/n$. It can be shown that $\mu(x)$ is a norm on $\mathbb{R}^2$ when $F(0)<1/2$, and the shape $A$ is the unit ball in this norm, that is, $A=\{x\in\mathbb{R}^2:\mu(x)\leq 1\}$. Clearly $A$ inherits all symmetries of $\mathbb{Z}^2$. Furthermore, the fact that $A$ is convex (which follows from subadditivity) implies that its boundary must fall somewhere between the diamond $\{x=(x_1,x_2)\in\mathbb{R}^2: |x_1|+|x_2|=\mu^{-1}\}$ and the cube with side $2\mu^{-1}$ centered at the origin. Constant passage times give the diamond as asymptotic shape, but characterizing $A$ for non-trivial distributions has turned out to be very difficult. Early simulations for exponentially distributed passage times indicated that $A$ was a circle; see e.g.\ \cite{Ed61, R73}. Kesten \cite{Kes86} however shows that, for a certain class of distributions including the exponential, $A$ is not a ball in sufficiently high dimension  (where ''sufficiently high`` means $d\geq 10^6$ for the exponential distribution), and this of course speaks against a circle in two dimensions. Furthermore, due to the underlying lattice structure there is no theoretical basis for a rotationally invariant shape.

In \cite{DL81} it is shown that $A$ has a flat edge (that is, its boundary contains a straight line segment) when $F$ is a shifted geometric distribution with $\PP(\tau(e)=k)=p(1-p)^{k-1}$ for $k\geq 1$ and $p$ sufficiently close to 1. Marchand \cite{M02} generalizes this by showing that any distribution $F$ with inf supp$(F)=r>0$ and $F(r)\geq p^{\sss \rm{dir}}_c$ gives rise to a flat edge, where $p_c^{\sss \rm{dir}}$ denotes the critical value for directed bond percolation on $\mathbb{Z}^2$. It is believed that $A$ is strictly convex (and hence does not have flat edges) as soon as
\begin{equation}\label{eq:F_conv_cond}
r=0 \mbox{ or } F(r)<p^{\sss \rm{dir}}_c,
\end{equation}
but so far strict convexity has not been proved for any distribution. In \cite{DamHoc11} however, distributions (with atoms) are constructed that give non-polygonal shapes, that is, shapes where the boundary does not consist only of flat edges.

As for comparing shapes for different distributions, the fact that $\mu^{\sss \widetilde{F}}<\mu^{\sss F}$ when $\widetilde{F}$ is more variable than $F$ can be generalized to any direction for distributions with inf supp$(F)$=0. This means that the asymptotic shape for $F$ is then strictly smaller than for $\widetilde{F}$. When inf supp$(F)=r>0$, the shape contains flat edges if $F(r)$ is large, as mentioned above. If in this case also inf supp$(\widetilde{F})=r$, then the shapes coincide along the flat edges, but otherwise and in other directions the shape for $\widetilde{F}$ is strictly larger than the shape for $F$; see \cite{M02} for details.

We will investigate the growing shape from a number of different aspects: For constant passage times, the time constant $\mu_\theta$ increases with $\theta$ so that the shape is smallest along the diagonals. \emph{Is this true in general? How does the shape depend on properties of the passage time distribution? Can we quantify the deviations from a circle?}

The asymptotic shape of first passage percolation has received attention also in statistical physics. There the focus has been on the fluctuations of $\bar{B}(t)$ around its mean shape. First passage percolation is conjectured to belong to a class of growth models analyzed by Kardar, Parisi and Zhang in \cite{KPZ86}. This means that the fluctuations of the interface between $\bar{B}(t)$ and its complement is believed to be of the order $t^{\chi}$ with $\chi=1/3$, so that hence Var$(T(\mathbf{0}, \mathbf{n}))$ is of the order $n^{2/3}$. The analysis also involves an exponent $\xi$ such that the fluctuations of the location where $\bar{B}(t)$ hits a line at distance $s$ from the origin is of the order $s^{\xi}$. Heuristical arguments suggest that the exponents $\chi$ and $\xi$ are related in that $\chi=2\xi-1$, which implies that $\xi=2/3$. Both exponents are believed to be independent of the direction and of the passage time distribution $F$. Our simulations give support for this, and the predicted values of the exponents (see also references below to earlier simulations).

Rigorous results related to the above conjectures concern Var$(T(\mathbf{0}, \mathbf{n}))$ but are quite weak. In \cite{Kes93}, it is shown that, if $F$ has finite second moment, then $c_1\leq \mbox{Var}(T(\mathbf{0}, \mathbf{n}))\leq c_2 n$. As for upper bounds, the best bound to date is  Var$(T(\mathbf{0}, \mathbf{n}))=O(n/\log n)$, which is proved for uniform distributions in \cite{BKS03} and for a larger class of (non-atomic) distributions in \cite{BR06}. As for lower bounds, the best result is that Var$(T(\mathbf{0}, \mathbf{n}))$ is at least of the order $\log n$ for distributions $F$ that satisfy (\ref{eq:F_conv_cond}); see \cite{NP95}.

Previous simulation studies on first passage percolation are limited, both in number and in scope. The early references \cite{Ed61, R73} have already been mentioned. There the shape is simulated for exponentially distributed passage times and is seen to be reminiscent of a circle. These simulations however were substantially restricted by limitations in computer power. In the theoretical physics literature, the so called Eden growth model \cite{Ed61} has received a lot of attention and one version of the model is equivalent to first passage percolation with exponential passage times. The interest mainly revolves around the scaling properties of the boundary of the growing infected set, see e.g.\ \cite{KS91}, and the model has also been subject to simulation studies from this aspect. The results confirm the conjecture for the exponents $\chi$ and $\xi$ described above; see e.g.\ \cite{ZS86}. The simulations are not performed on the version of the model that is equivalent to first passage percolation, but the different versions are closely related and all belong to the universality class defined in \cite{KPZ86}, so the conclusions are still relevant in this context.

\section{Simulations}

Our simulations are designed so that, in each run, we keep track of the infected set until all lines at distance $x_{\max}$ from the origin have been hit. Here $x_{\sss max}$ varies between 10.000 and 20.000 depending on the distribution of the passage times, see Table \ref{tab_summary}. We record all hitting times to lines in five symmetry directions, including the axis and the diagonal, and, for each such direction, we calculate average hitting times over all symmetries (the averages are hence based on four values in the axis direction and along the diagonal, and on eight values for the other directions). The directions in between the axis and the diagonal are, in the first octant, defined by the lines $4x+y=c$, $2x+y=c$ and $4x+3y=c$ for increasing values of $c$, and the corresponding directional angles with the $x$-axis are approximately 14, 27 and 37 degrees, respectively. The lines run through integer points, which makes it slightly easier to record the hitting times. We also keep track of the hitting points on the lines. The simulations are performed on a cluster of 28 Linux machines, each with 12 kernels, running at 2 GHz.

We study several different continuous passage time distributions. Note that all distributions are scalable in the sense that multiplication by a constant (which allows us to tune e.g.\ the mean) gives rise to a distribution in the same distribution class. 

\begin{itemize}
 \item[1.] The exponential distribution with mean 1. Because of the memoryless property of the exponential distribution, the simulations here are much faster than for other distributions and are therefore more extensive. For comparison, we mention that one simulation with $x_{\sss max}=20.000$ takes 22 minutes for the exponential distribution and 37 hours for the uniform distribution on [0,1] -- the times hence differ roughly by a factor 100.
 \item[2.] Gamma distributions, $\Gamma(k,k)$, with $k=2,3,4$. The mean value is 1 for all $k$, but the variance decreases as $k$ increases. Since these Gamma distributions are convolutions of exponential distributions, we can exploit the memoryless property also here to speed up the simulations.
 \item[3.] Uniform distributions on $[0,1]$ and on $[0.1,0.9]$. 
 \item[4.] Power-law distributions with support on $(0,\infty)$ and tail exponent $\gamma=1,2,3,4$. More specifically, the passage times are sampled from a scaled Fisher distribution defined by $\widetilde{\rm{Fi}}(\gamma)=\rm{Fi}(2,2\gamma)/\gamma$. We then have $1-F(x)=(1+x)^{-\gamma}$. 
\end{itemize}

Let $\tau,\tau_1,\ldots, \tau_4$ be i.i.d.\ with distribution $F$ and define $\E_4[\tau]=\E[\min\{\tau_1,...,\tau_4\}]$ and $\E_4[\tau^2]=\E[\min\{\tau_1^2,...,\tau_4^2\}]$. We recall from Section 1 that $\E_4[\tau]<\infty$ guarantees the existence of a time constant and $\E_4[\tau^2]<\infty$ the existence of an asymptotic shape. Table \ref{tab_summary} contains the values of these quantities along with the mean and the variance of the distributions and a summary of the scope of the performed simulations. Based on the findings for these large-scale simulations, we have also performed additional smaller simulations for translated versions of the distributions; see below for further details.


\begin{table}
\centering
\begin{tabular}{|l | c | r | l | l | l | l | }
\hhline{~|------}
\multicolumn{1}{c|}{} & $x_{\sss max}$ & $\#$ runs & $\E[\tau]$ & Var($\tau$) & $\E_4[\tau]$ & $\E_4[\tau^2]$ \\ \hline
Exp(1) & 20.000 & 20.000 & 1 & 1 & 0.25 & 0.125 \\
$\Gamma(2,2)$ & 20.000 & 5000 & 1 & 0.5 & 0.4023 & 0.2256\\
$\Gamma(3,3)$ & 20.000 & 5000 & 1 & 0.333 & 0.4887 & 0.2957\\
$\Gamma(4,4)$ & 18.000 & 5000 & 1 & 0.25 & 0.5457 & 0.3480\\
U(0,1) & 20.000 & 500 & 0.5 & 0.0833 & 0.2 & 0.0667\\
U(0.1,0.9) & 20.000 & 500 & 0.5 & 0.0533 & 0.26 & 0.0847\\
$\widetilde{\rm{Fi}}(1)$ & 10.000 & 1000 & - & - & 0.3333 & 0.3333\\
$\widetilde{\rm{Fi}}(2)$ & 10.000 & 1000 & 1 & - & 0.1429 & 0.0476\\
$\widetilde{\rm{Fi}}(3)$ & 10.000 & 1000 & 0.5 & 0.75 & 0.0909 & 0.0182\\
$\widetilde{\rm{Fi}}(4)$ & 10.000 & 1000 & 0.3333 & 0.2222 & 0.0667 & 0.0095\\ \hline
\end{tabular}
\caption{Summary of simulated passage time distributions.}
\label{tab_summary}
\end{table}

\subsection{Time constants}\label{sim_tc}

Here we first describe the estimation procedure for the time constants and then the results of the simulations.

\subsubsection*{Estimation}

We first estimate the time constants in the directions described above. The simplest estimation procedure would be to take the scaled average passage times at distance $x_{\sss max}$ from the origin. However, in order to increase the accuracy, we extrapolate to infinity as follows: Take $\theta=0$ (the procedure is the same in all directions, so we describe it for the axis direction) and let $\hat{\mu}_n$ denote the average passage time to the line $x=n$ divided by $n$. To estimate $\mu$, we use a weighted regression model with
$$
\E[\hat{\mu}_n]=\mu+cn^{-2/3},
$$
where the weight of each $\hat{\mu}_n$ is given by the inverse standard deviation. We pick the sample points at distance 500 starting at $n=1000$ (since the observation at $n=500$ is seen to be an outlier). In view of the predictions for the variance of the passage times, the factor $n^{-2/3}$ is a natural choice. It also produces a very good fit with data, as can be seen from Figure \ref{fig:reg}, which shows a fitted regression line for the exponential distribution (the other distributions give rise to similar pictures). The coefficients of determination for our regressions are throughout very high, with values ranging between 0.9974 and 1. In Table \ref{tab_est}, we summarize the obtained estimates.

To calculate the standard errors of the estimates, we note that the values at different observation points in the regression are positively correlated and the resulting estimate of the standard deviation is therefore likely to be too small. To get more conservative values, we proceed as follows: For each distribution, the data is stored in a number of files, each one containing average passage times over a number of runs (the number of files and the number of runs in each file varies between the distributions). We do one extrapolation per file, as described above, and then calculate the standard deviation in the resulting set of estimates. Comparing with the standard errors based on one single extrapolation, the obtained values are roughly a factor 10 larger. As this second method does not use all information optimally, the true values most likely lie somewhere in between the obtained values. Table \ref{tab_est} contains the conservative estimates based on repeated extrapolations. We could of course also base our estimates of $\mu$ on the repeated extrapolations by calculating the average estimate, but this makes little difference -- the resulting estimates generally coincide up to the sixth decimal.

\subsubsection*{Results} 

\begin{figure}[h]
\centering
\begin{minipage}{.48\textwidth}
  \begin{center}
  \includegraphics[width=7cm]{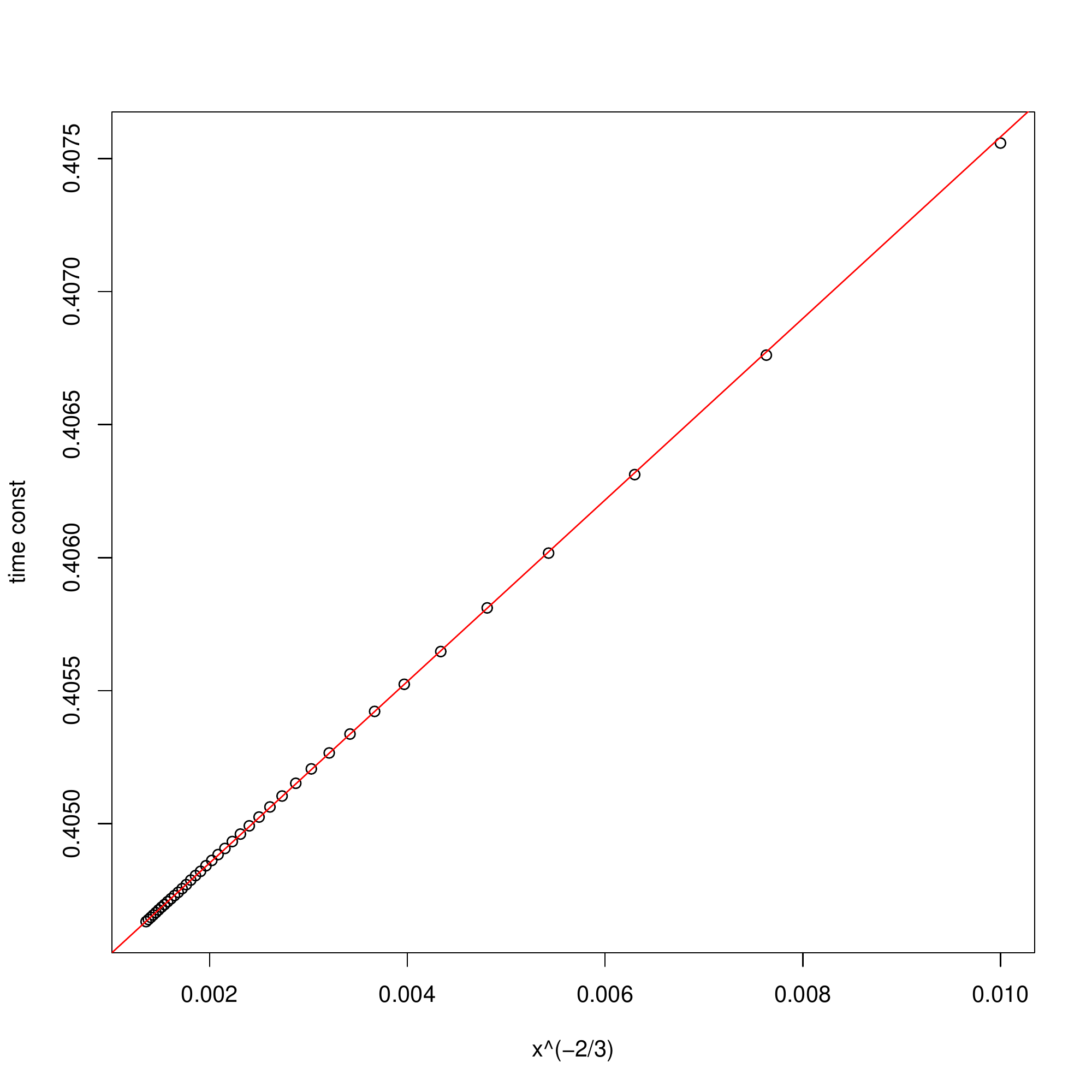}
  \end{center}
  \caption{Regression line for estimation of the time constant for the exponential distribution.}
  \label{fig:reg}
\end{minipage}
\hspace{0.2cm}
\begin{minipage}{.48\textwidth}
  \begin{center}
  \includegraphics[width=8cm]{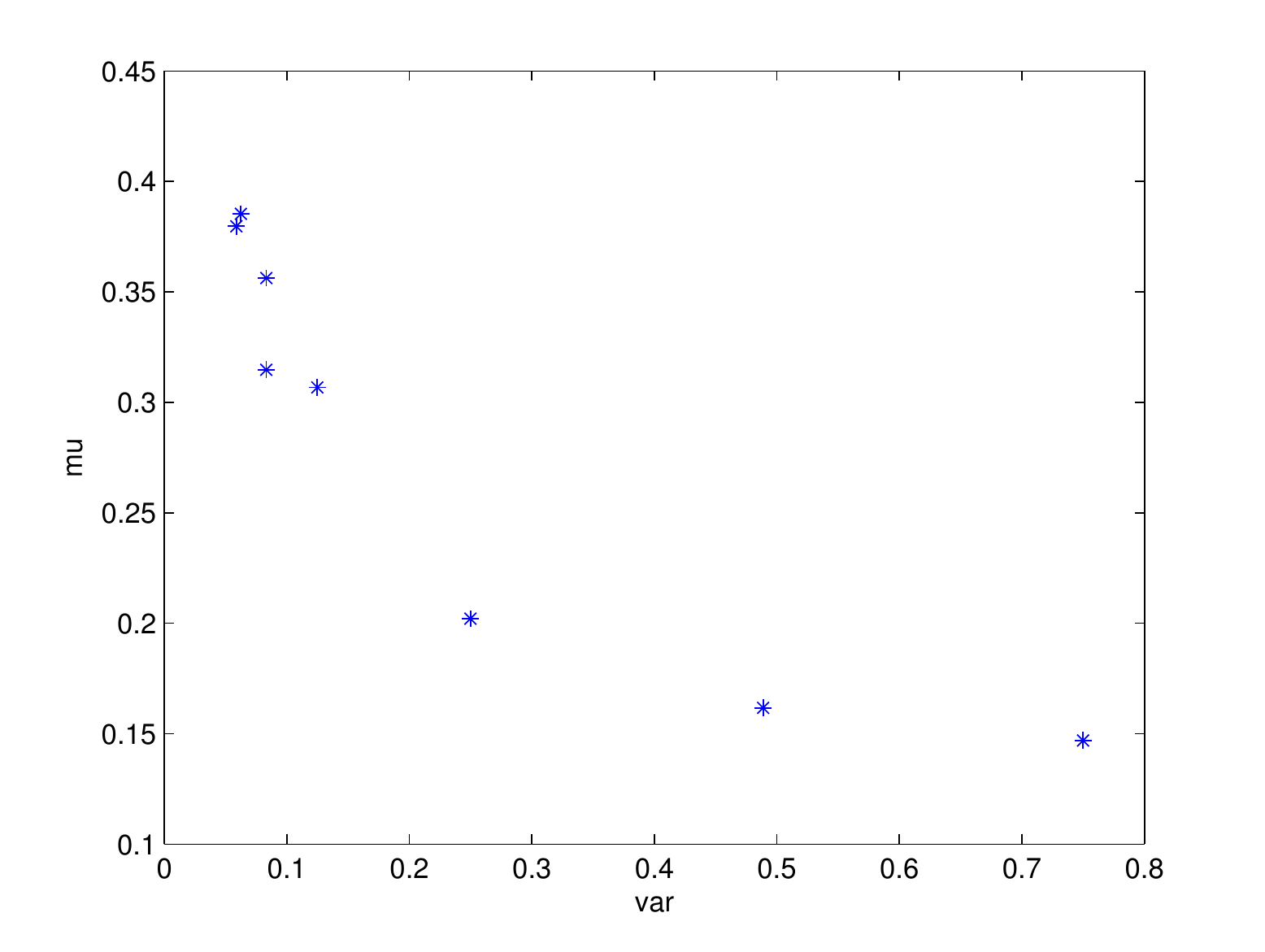}
  \end{center}
  \caption{Estimated time constant plotted against the variance for scaled distributions with mean 0.5 and finite variance.}
  \label{fig:tc_var}
\end{minipage}
\end{figure}

We first study the time constant in a fixed direction, e.g.\ the axis direction. Can we find a quantity that determines its value? As pointed out in Section 2.1, two distributions with the same mean does not necessarily give rise to the same time constant and, scaling the simulated distributions so that they all have mean 0.5, we see that the corresponding values of $\mu=\mu_{\sss 0}$ are indeed very different. Does the variance play a role? Figure \ref{fig:tc_var} shows a plot of the time constants against the variance for the scaled distributions (with mean 0.5) and there does seem to be a tendency that a larger variance implies a smaller time constant (i.e.\ faster growth of the infection), but the relation is far from linear and does not explain all variation in the time constant. Comparing e.g.\ U(0,1) and $\Gamma(3,3)/2$, that both have mean 0.5 and variance 1/12, we see that the corresponding time constants are quite different -- 0.3131 and 0.3563, respectively. Hence also the mean and the variance together are not sufficient for determining the value of $\mu$. 

\begin{figure}[h]
\centering \mbox{\subfigure[$\theta=0$]{\includegraphics[height=5.7cm]{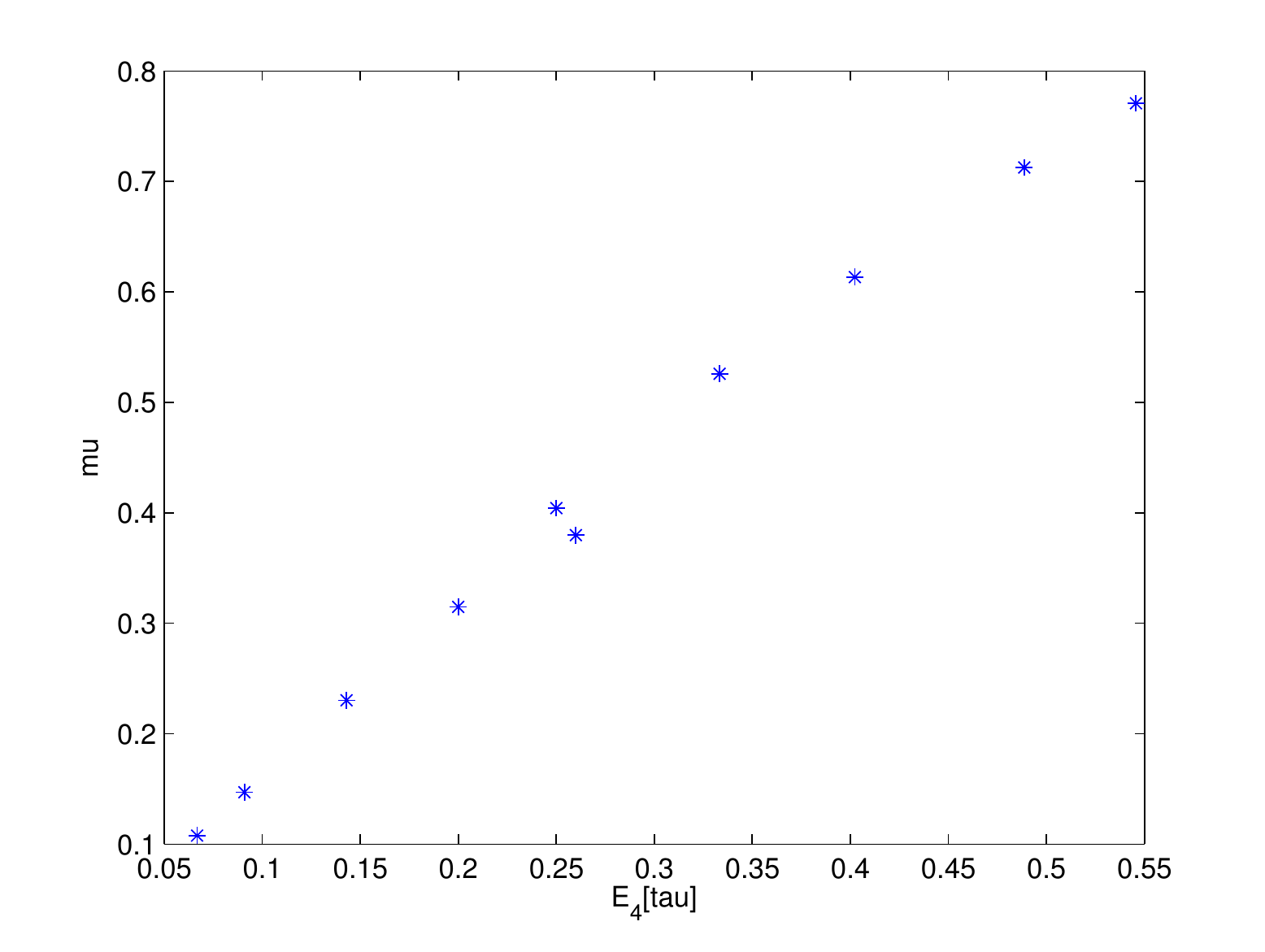}}}\quad
\mbox{\subfigure[$\theta=45$]{\includegraphics[height=5.7cm]{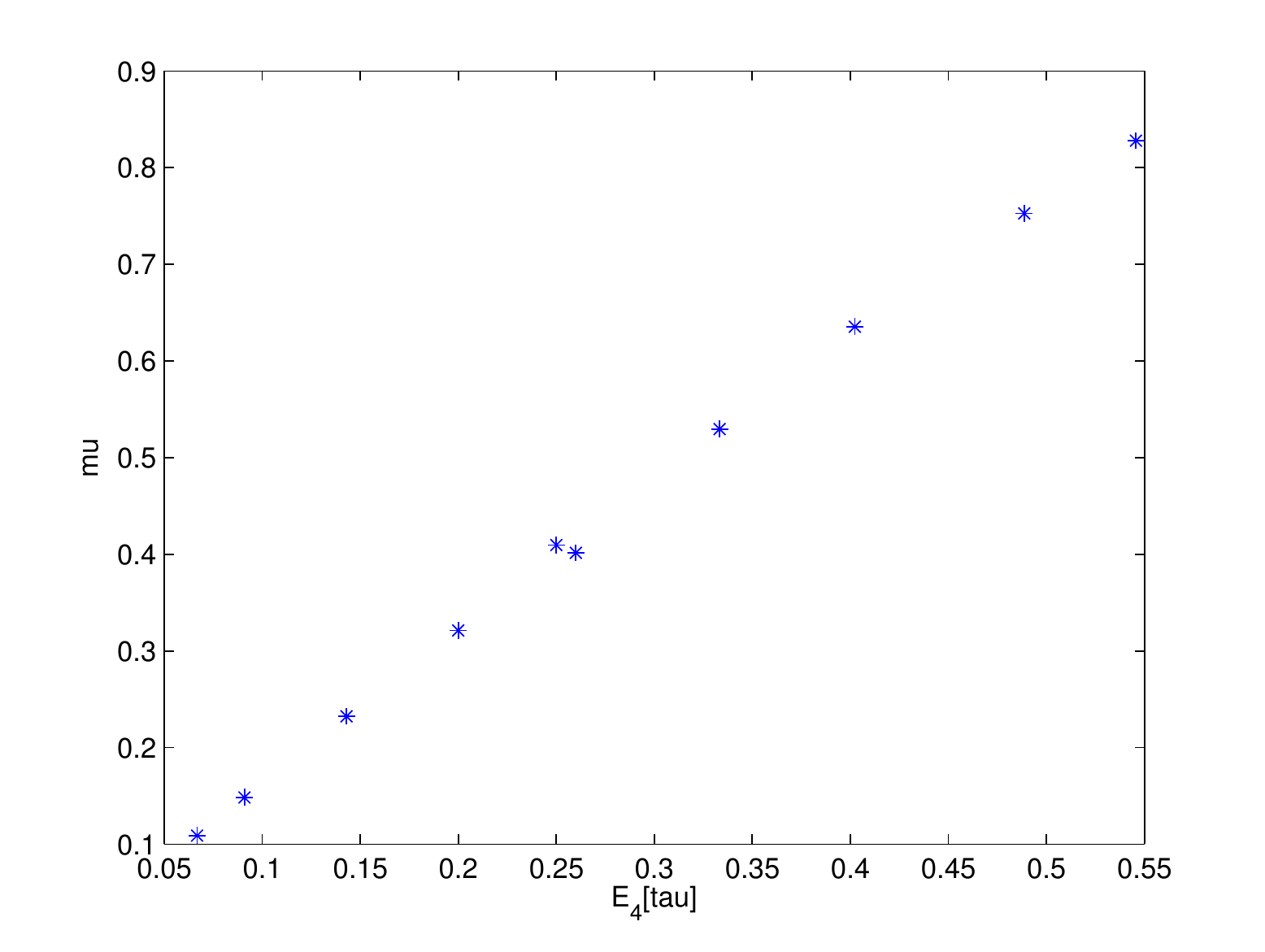}}}\caption{Estimated time constants plotted against $\E_4[\tau]$ for the distributions in Table \ref{tab_summary}.}\label{fig:tc_E4}
\end{figure}

We next turn to $\E_4[\tau]$, which is a natural quantity to look at since finiteness of this expectation guarantees the existence of a time constant. For U(0,1) and $\Gamma(3,3)/2$, this equals 0.2 and 0.2444, respectively, which indeed offers a possible explanation of the difference between the time constants. In Figure \ref{fig:tc_E4}, the time constants along the axis and along the diagonal are plotted against $\E_4[\tau]$ for all simulated distributions (the distibutions are not scaled to have the same mean) and in both directions we observe an almost perfect linear relation. Pictures in the other directions look similar. We note that the observation for U(0.1,0.9) slightly deviates from the linear pattern. Removing this observation and fitting regression lines through the origin to the other observations gives $\mu_0=1.5\cdot \E_4[\tau]$ and $\mu_{45}=1.56\cdot \E_4[\tau]$, respectively, with increasing values between 1.5 and 1.56 for the slopes at the intermediate angles. 

\begin{table}[t]
\centering
\begin{tabular}{|l | l | l | l | l | l | l|}
\hhline{~~|-----} \multicolumn{2}{c|}{} & $\theta=0$ & $\theta=14$ & $\theta=27$ & $\theta=37$ & $\theta=45$\\ \hline
\multirow{2}{*}{Exp(1)} & $\mu_{\sss\theta}$ & 0.4041702 & 0.4052820 & 0.4074979 & 0.4090543 & 0.4095065\\
 & se & 0.0000027 & 0.0000018 &  0.0000017 &  0.0000019 &  0.0000021\\ \hline
\multirow{2}{*}{$\Gamma$(2,2)} & $\mu_{\sss\theta}$ & 0.6134115 & 0.6176678 & 0.6265413 & 0.6332689 & 0.6352802 \\
 & se & 0.0000034 & 0.0000039 &  0.0000038 &  0.0000036 &  0.0000051\\ \hline
\multirow{2}{*}{$\Gamma$(3,3)} & $\mu_{\sss\theta}$ & 0.7125932 & 0.7198929 & 0.7358333 & 0.7487188 & 0.7527526 \\
 & se & 0.0000055 & 0.0000037 &  0.0000027 &  0.0000029 &  0.0000034\\ \hline
\multirow{2}{*}{$\Gamma$(4,4)} & $\mu_{\sss\theta}$ & 0.7706122 & 0.7804740 & 0.8026975 & 0.8216638 & 0.8278508 \\
 & se & 0.0000050 & 0.0000036 &  0.0000041 &  0.0000034 &  0.0000045\\ \hline
\multirow{2}{*}{U(0,1)} & $\mu_{\sss\theta}$ & 0.3131102 & 0.3147531 & 0.3180768 & 0.3204898 & 0.3212093 \\
 & se & 0.0000109 & 0.0000079 &  0.0000068 &  0.0000087 &  0.0000010\\ \hline
\multirow{2}{*}{U(0.1,0.9)} & $\mu_{\sss\theta}$ & 0.3796196 & 0.3836337 & 0.3923484 & 0.3993590 & 0.4015409 \\
 & se & 0.0000087 & 0.0000047 & 0.0000056 &  0.0000058 &  0.0000076\\ \hline
\multirow{2}{*}{$\widetilde{\rm{Fi}}(1)$} & $\mu_{\sss\theta}$ & 0.5257341 & 0.5265118 & 0.5281133 & 0.5291650 & 0.5294726 \\
 & se & 0.0000286 & 0.0000207 &  0.0000205 &  0.0000192 &  0.0000266\\ \hline
\multirow{2}{*}{$\widetilde{\rm{Fi}}(2)$} & $\mu_{\sss\theta}$ & 0.2301307 & 0.2305824 & 0.2314898 & 0.2321377 & 0.2323273 \\
 & se & 0.0000106 & 0.0000080 &  0.0000081 & 0.0000077 &  0.0000116\\ \hline
\multirow{2}{*}{$\widetilde{\rm{Fi}}(3)$} & $\mu_{\sss\theta}$ & 0.1468621 & 0.1471894 & 0.1478335 & 0.1482962 & 0.1484319 \\
 & se & 0.0000070 & 0.0000046 &  0.0000049 &  0.0000054 &  0.0000073\\ \hline
\multirow{2}{*}{$\widetilde{\rm{Fi}}(4)$} & $\mu_{\sss\theta}$ & 0.1077790 & 0.1080357 & 0.1085429 & 0.1088943 & 0.1089969 \\
 & se & 0.0000062 & 0.0000033 &  0.0000036 &  0.0000036 &  0.0000048\\ \hline
\end{tabular}
\caption{Estimates of directional time constants and corresponding conservative standard errors for the simulated distributions. Estimates for scaled distributions are obtained by scaling the estimates analogously.}
\label{tab_est}
\end{table}

\begin{figure}
\centering \mbox{\subfigure[$\theta=0$]{\includegraphics[height=5cm]{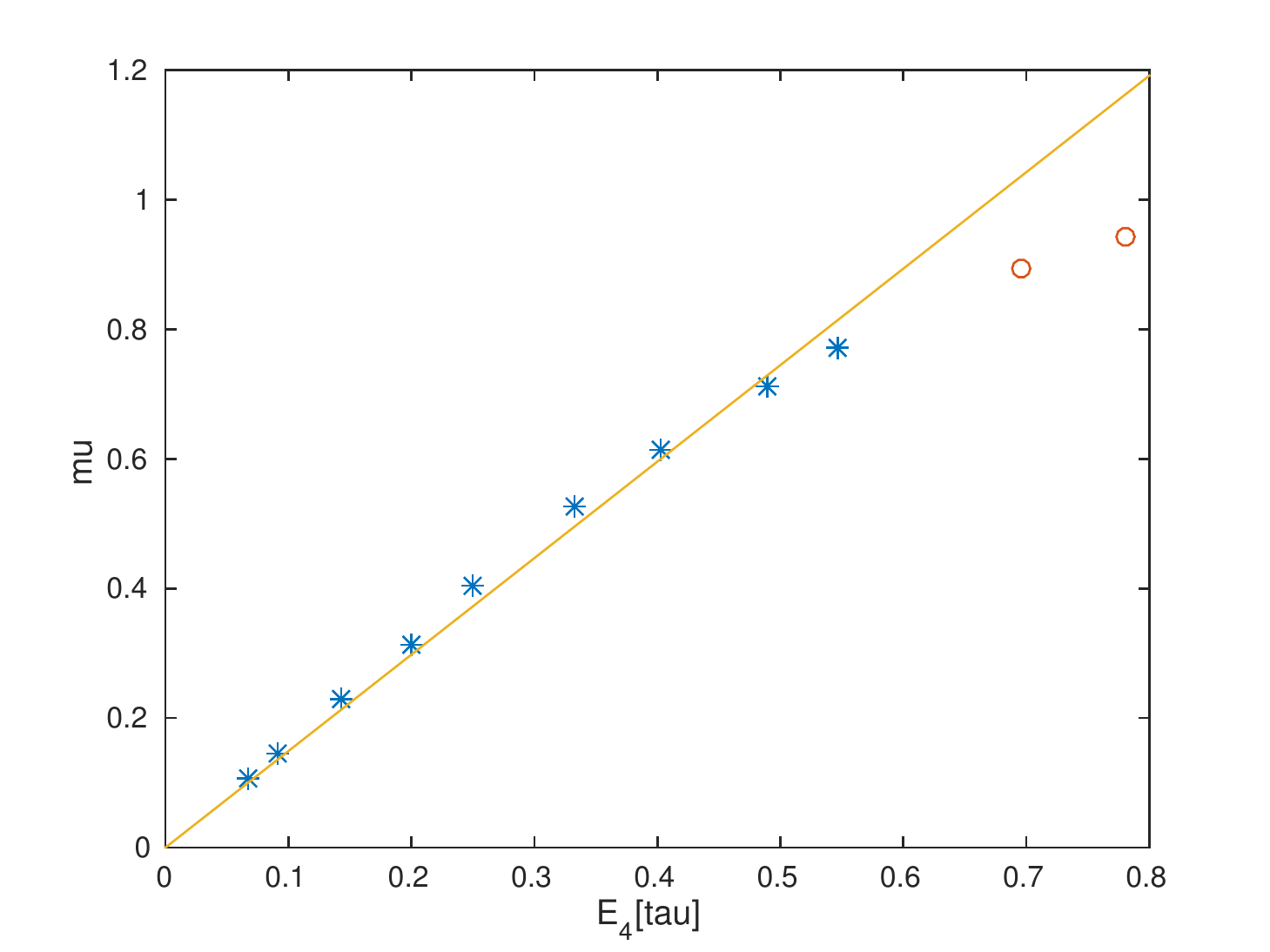}}}\quad
\mbox{\subfigure[$\theta=45$]{\includegraphics[height=5cm]{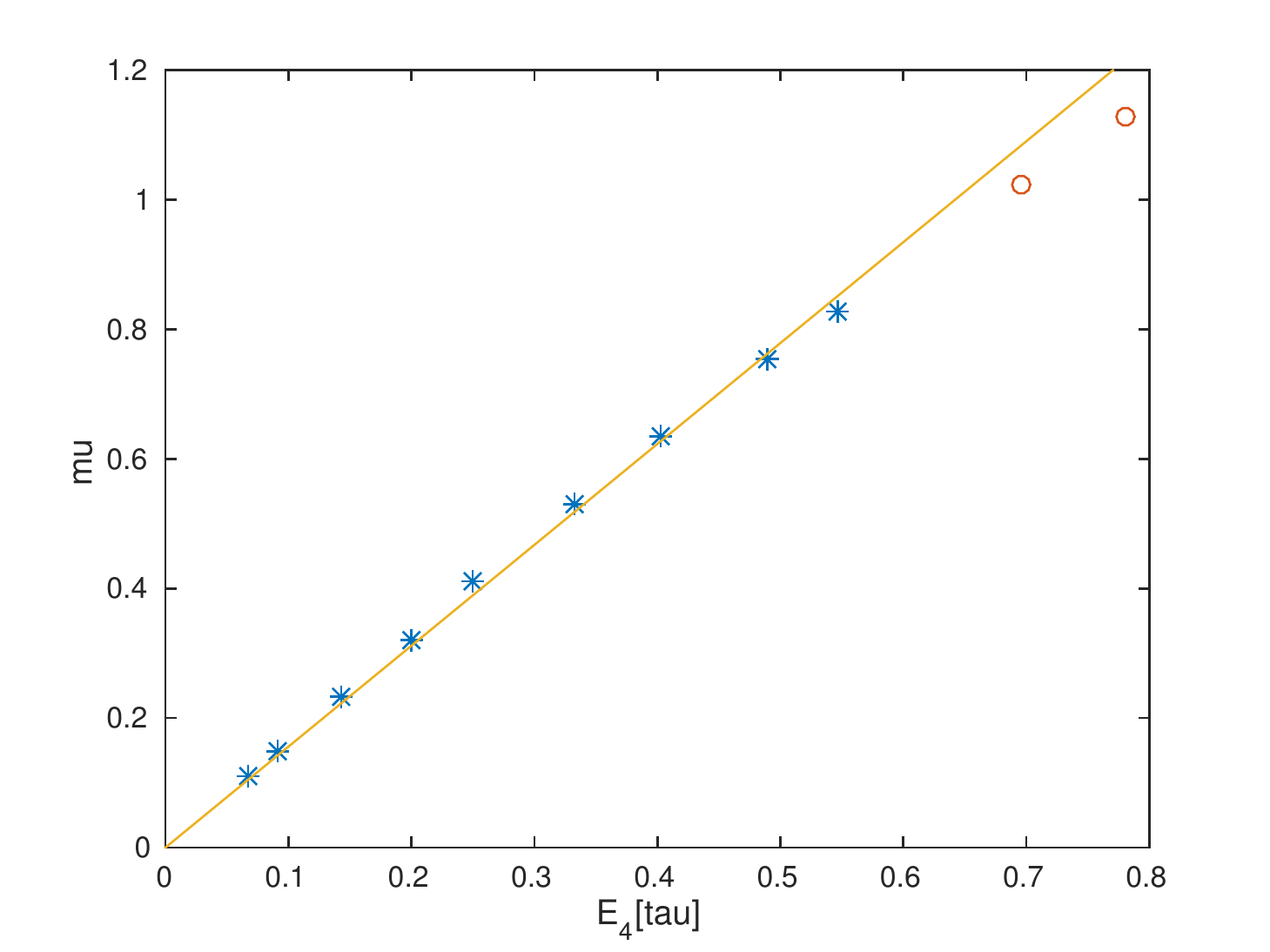}}}\caption{Estimated time constants plotted against $\E_4[\tau]$ for the distributions in Table \ref{tab_summary}, except U(0.1,0.9), (stars) and $\Gamma(k,k)$ for $k=10,20$ (circles). Regression lines fitted to the distributions in Table \ref{tab_summary}, except U(0.1,0.9).}\label{fig:tc_E4_spik}
\end{figure}

Our findings indicate that the time constant is closely related to the quantity $\E_4[\tau]$. A Gamma distribution and a Fisher distribution have very different properties and therefore it seems reasonable to believe that the above approximate relation between the time constant and $\E_4[\tau]$ is valid for a quite large class of distributions, e.g.\ absolutely continuous distributions $F$ with inf supp$(F)=0$ that are not too concentrated. If the distribution is very concentrated, e.g.\ with mass $1-\varepsilon$ on the interval $[a-\varepsilon,a+\varepsilon]$ for some $a>0$ and the rest of the mass smeared out over the positive axis, then the time constant and $\E_4[\tau]$ will both be close to $a$ and their ratio close to 1 (rather than 1.5) so it is clear that the relation is not completely universal. Figure \ref{fig:tc_E4_spik} is analogous to Figure \ref{fig:tc_E4} (with U(0.1,0.9) excluded), but we have included estimates also for $\Gamma(10,10)$ and $\Gamma(20,20)$ based on smaller simulations (500 runs with $x_{max}=10.000$). These distributions have variance $1/10^2$ and $1/20^2$, respectively, and we note that, for $\theta=0$, their values indeed fall below a line fitted to the other distributions while, for $\theta=45$, they still comply quite well with the other distributions (see below for some further comments on this). As for theoretical bounds, we believe that $\mu\geq \E_4[\tau]$ is true in great generality -- perhaps as soon as $F(0)<1/2$ so that $\mu>0$, and most likely when $F(0)=0$ -- but we do not have a proof of this. 

\begin{figure}
\centering \mbox{\subfigure[$\theta=0$]{\includegraphics[height=5cm]{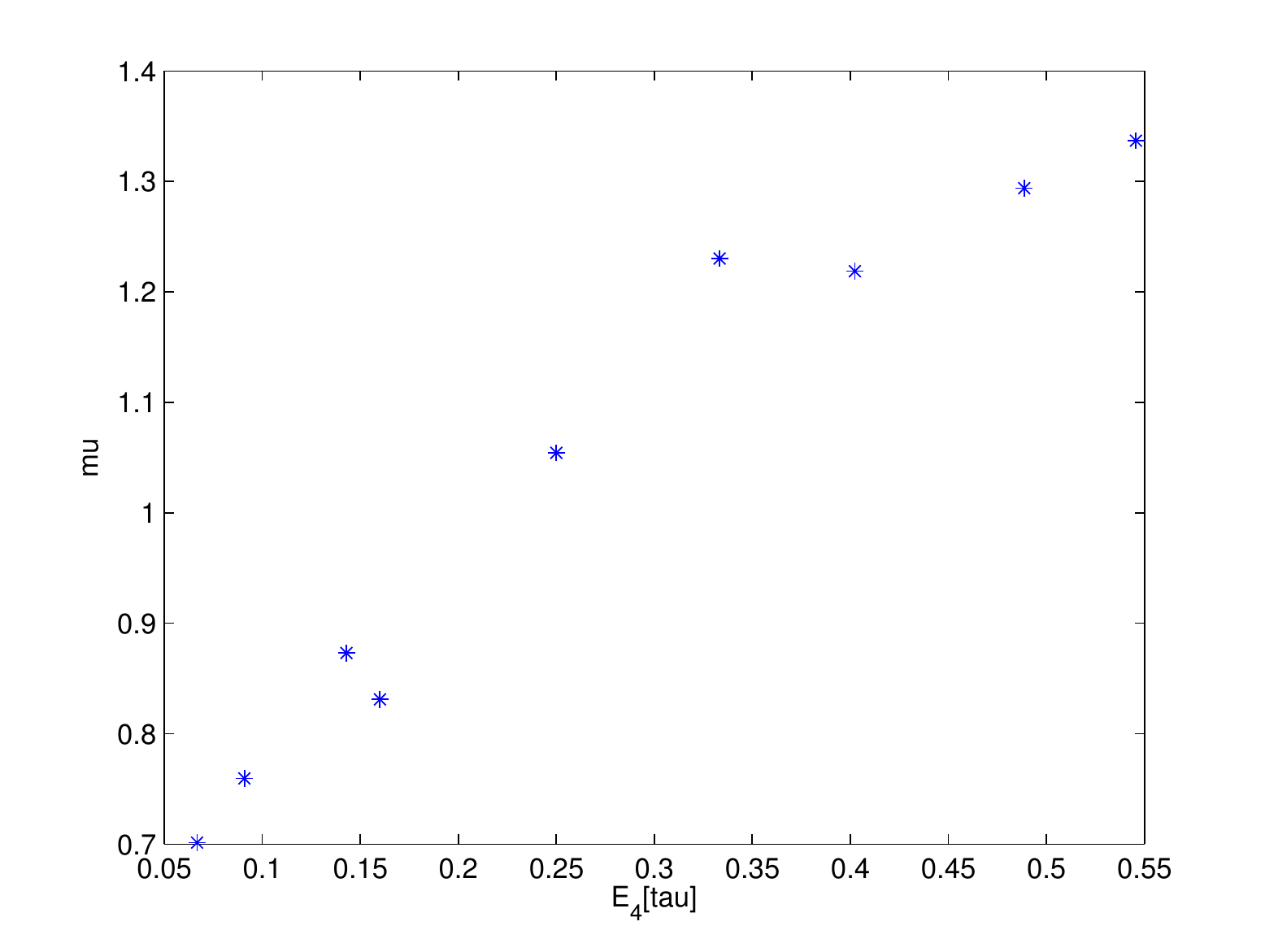}}}\quad
\mbox{\subfigure[$\theta=45$]{\includegraphics[height=5cm]{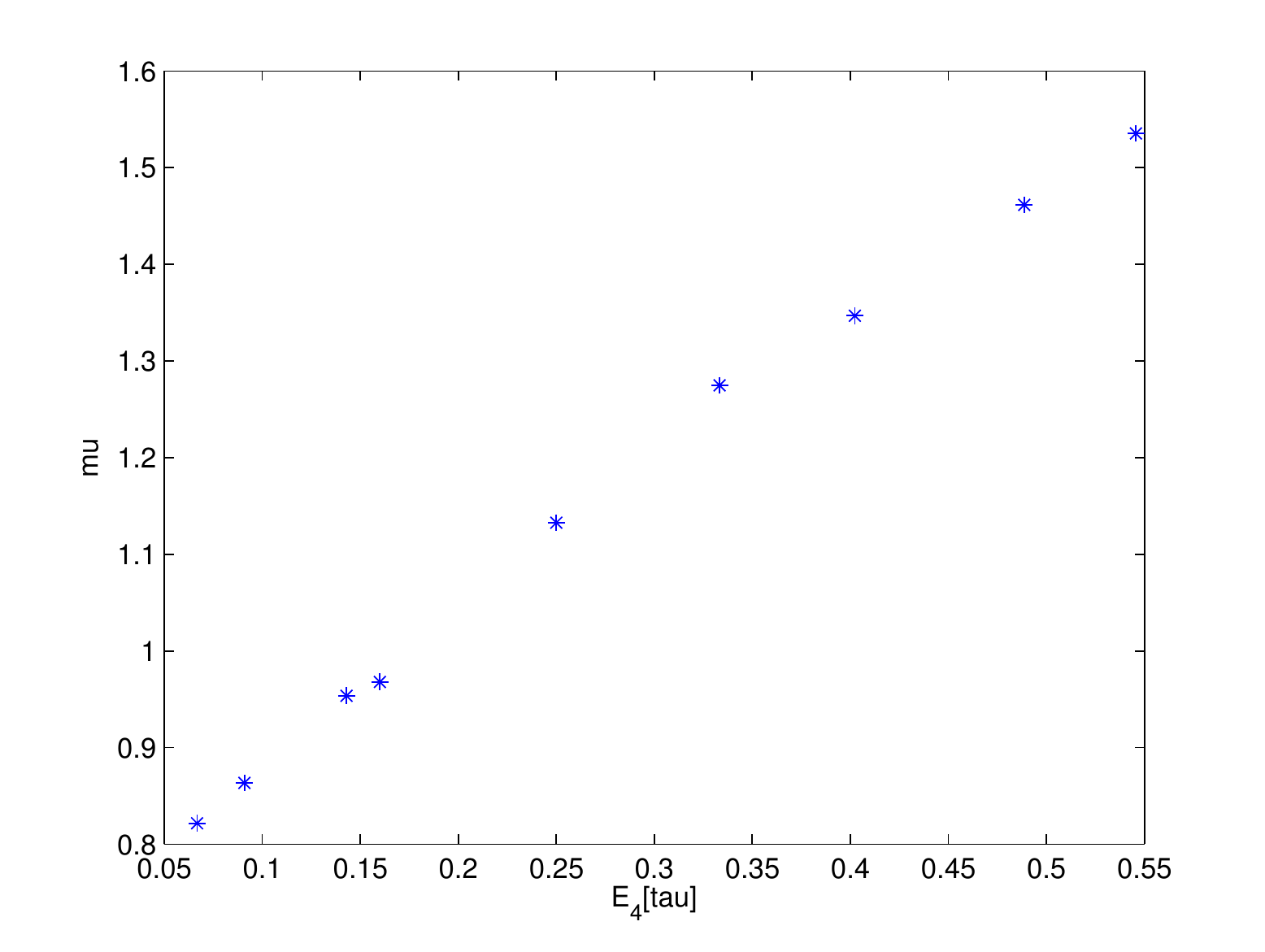}}}\caption{Estimated time constants for translated passage times $0.5+\tau$ plotted against $\E_4[\tau]$ for the distributions in Table \ref{tab_summary}, except U(0.1,0.9).}\label{fig:tc_E4_05}
\end{figure}

\begin{figure}
\centering \mbox{\subfigure[$\theta=0$]{\includegraphics[height=5cm]{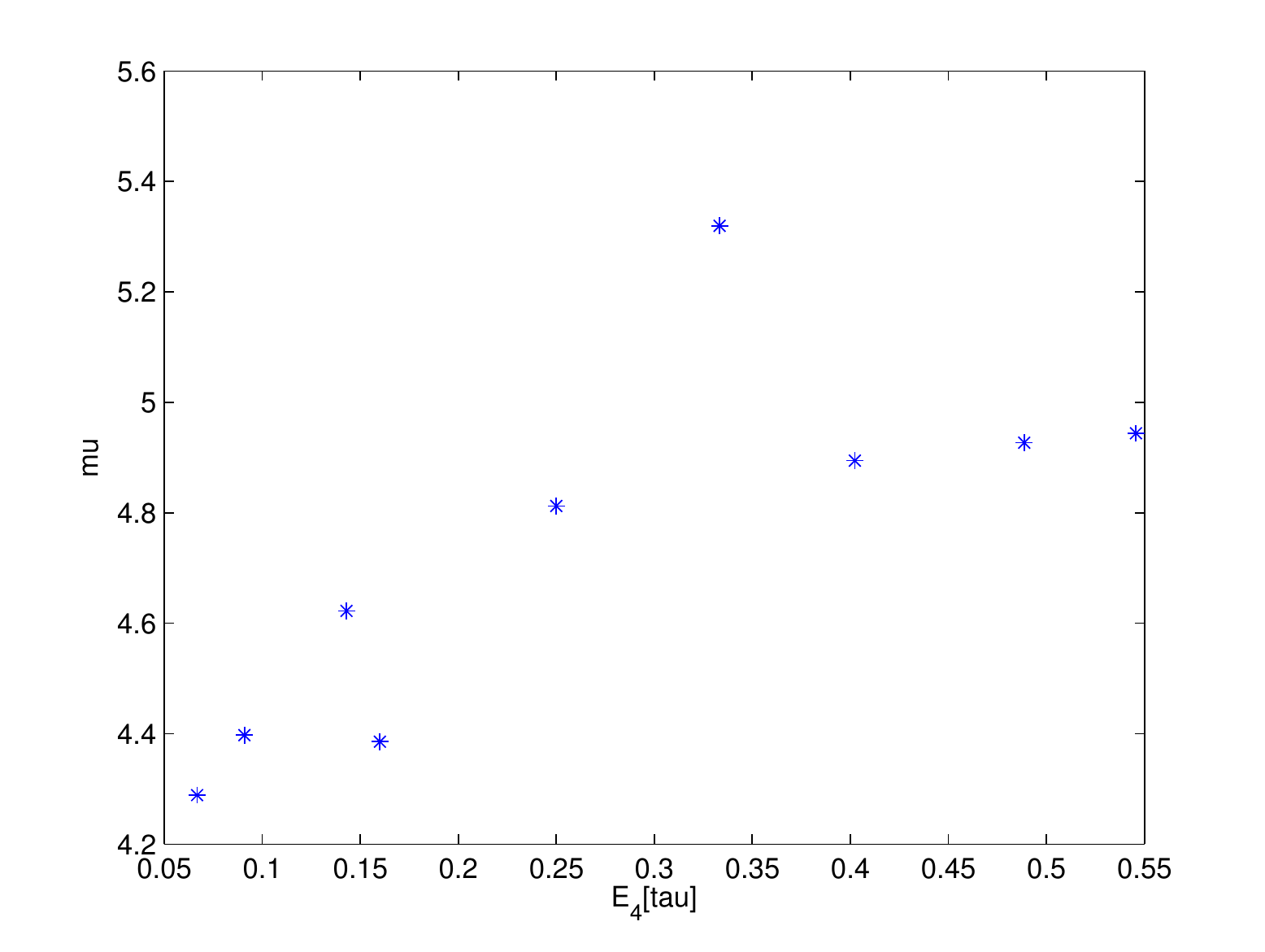}}}\quad
\mbox{\subfigure[$\theta=45$]{\includegraphics[height=5cm]{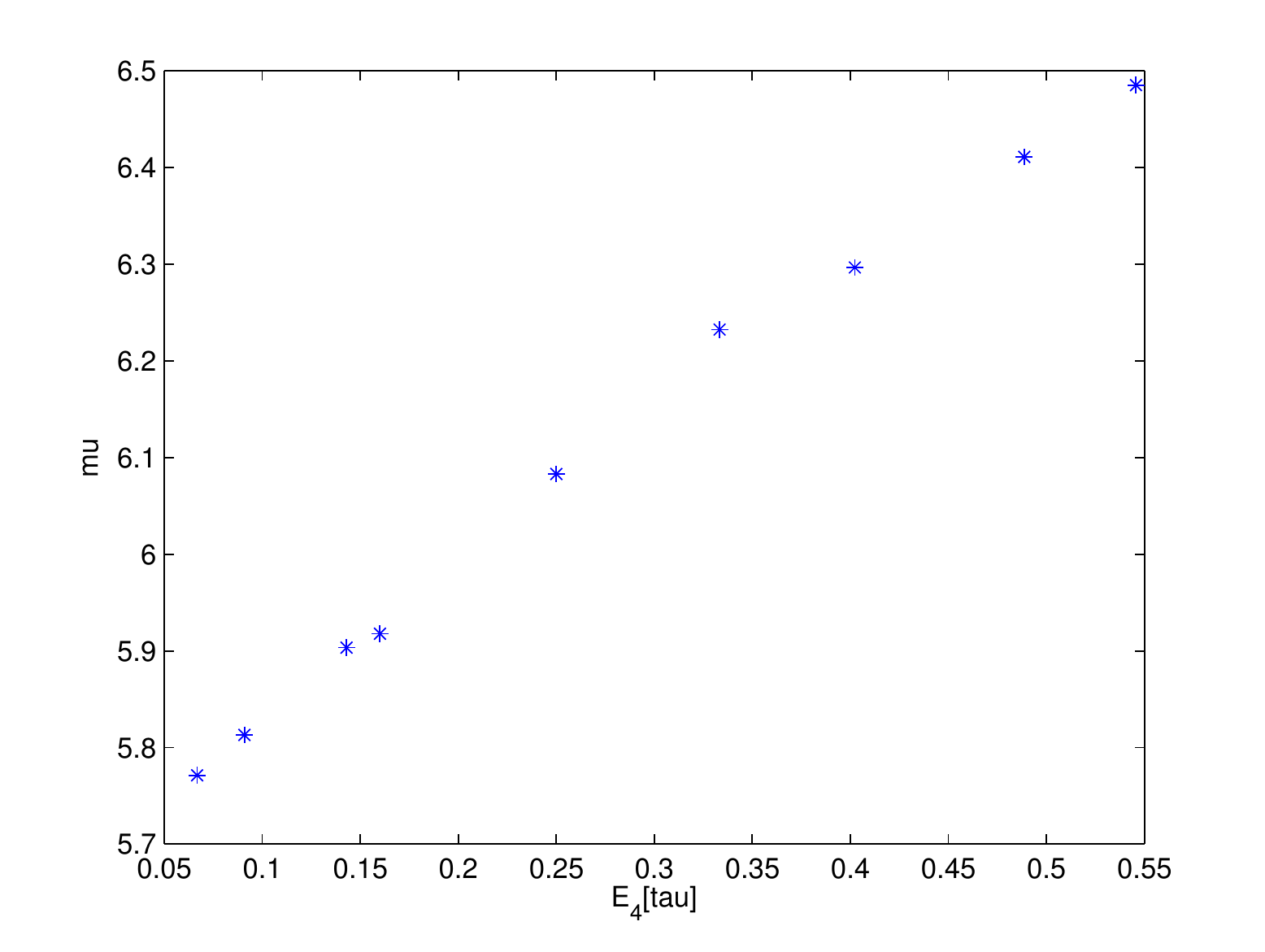}}}\caption{Estimated time constants for translated passage times $4+\tau$ plotted against $\E_4[\tau]$ for the distributions in Table \ref{tab_summary}, except U(0.1,0.9).}\label{fig:tc_E4_4}
\end{figure}

Figure \ref{fig:tc_E4} revealed that U(0.1,0.9) deviates slightly from the pattern for the distributions with inf supp$(F)=0$. What happens in general to the time constant when the support of the distribution is bounded away from 0? In order to investigate this, we have performed additional simulations for translated versions of our distributions, more specifically, we have simulated passage times distributed as $C+\tau$, where $C\in\{0.1,0.2,0.5,1,2,4\}$ and $\tau$ has one the the distributions in Table \ref{tab_summary} (excluding U(0.1,0.9)). These simulations are smaller and consist of 100 runs to $x_{max}=10.000$ for each translated distribution. Plotting the time constants against $\E_4[\tau]$ indicates that the relation becomes (i) less linear when $C$ increases for a given angle $\theta$ and (ii) more linear when $\theta$ increases for a given value of $C$. This is illustrated with pictures for $\theta=0$ and $\theta=45$ in Figure \ref{fig:tc_E4_05} ($C=0.5$) and Figure \ref{fig:tc_E4_4} ($C=4$). 

The fact that very concentrated distributions and translated distributions deviate from the pattern observed for more spread out distributions with inf supp$(F)=0$ may be explained in that, in both these cases, the model is regularized, meaning roughly that it contains less randomness relative to the average passage time. Thereby it comes closer to the situation with constant passage time. In both cases, there will be a tendency that the shortest paths are also the fastest, so that there is a stronger tendency for the infection to hit a point via the neighbor on the shortest path. Using $\E_4[\tau]$ as an explaining factor however is based on the assumption that the infection may invade a point via any of its four neighbors. The regularizing effect is most evident along the axis, where there is a unique shortest path, while it takes longer to set in along the diagonal, where there are more paths of optimal length. Quantifying what happens in the transition from the regime where $\E_4[\tau]$ explains the time constant to the case with constant passage times seems difficult.

How does the time constant depend on the direction? From Table \ref{tab_est} we see that, for each given distribution, the estimated time constant increases with $\theta$. Comparing two consecutive angles, the differences are of the order $10^{-3}$, while the standard deviations are of the order $10^{-6}$, and the differences are hence statistically significant. We conclude that the growth is fastest along the axis and slowest along the diagonal. The conclusion is the same for the additional simulations for translated distributions.

\subsection{Asymptotic shape}\label{sim_shape}

\begin{figure}
\centering \mbox{\subfigure[$k=1$]{\includegraphics[height=7.5cm]{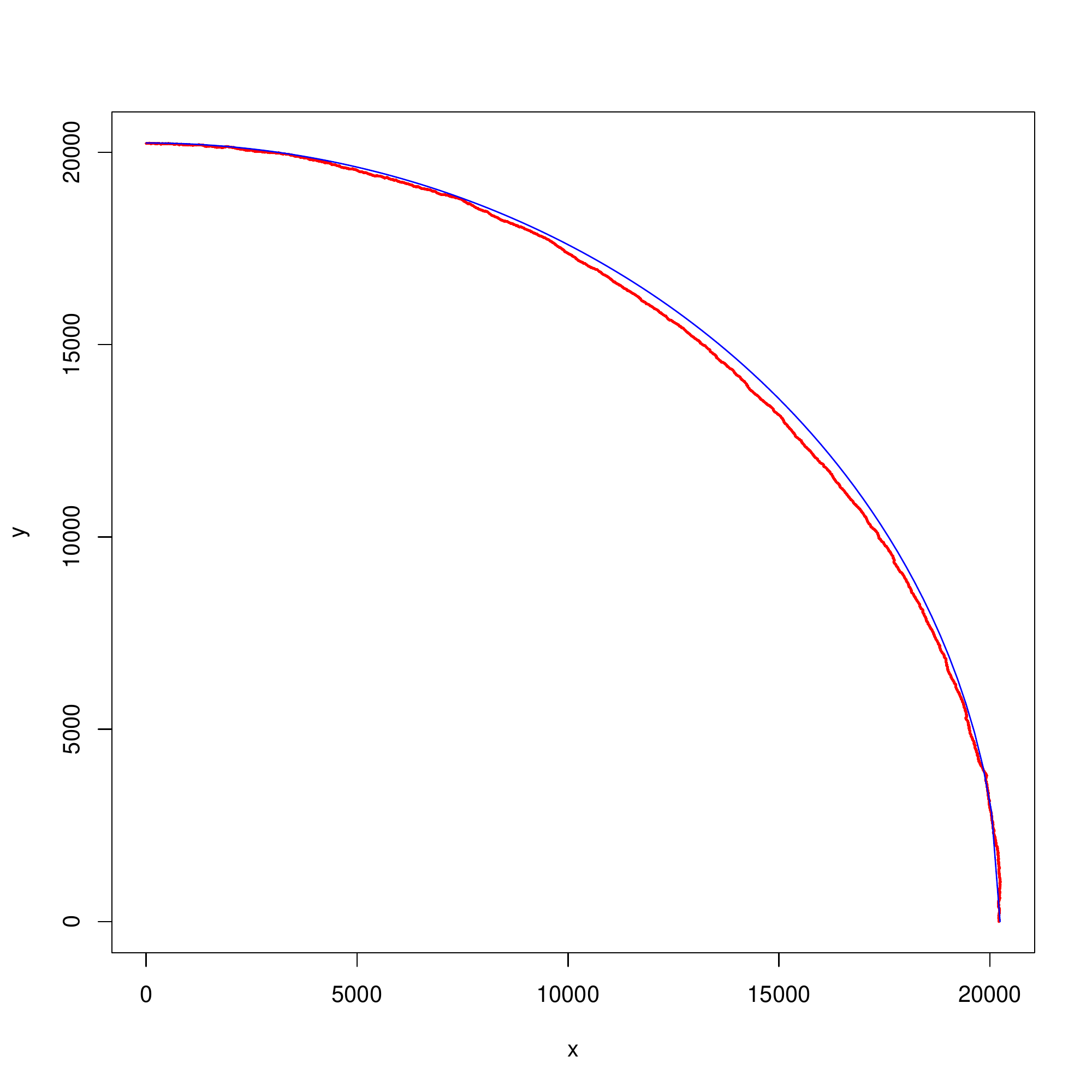}}}\quad
\mbox{\subfigure[$k=2$]{\includegraphics[height=7.5cm]{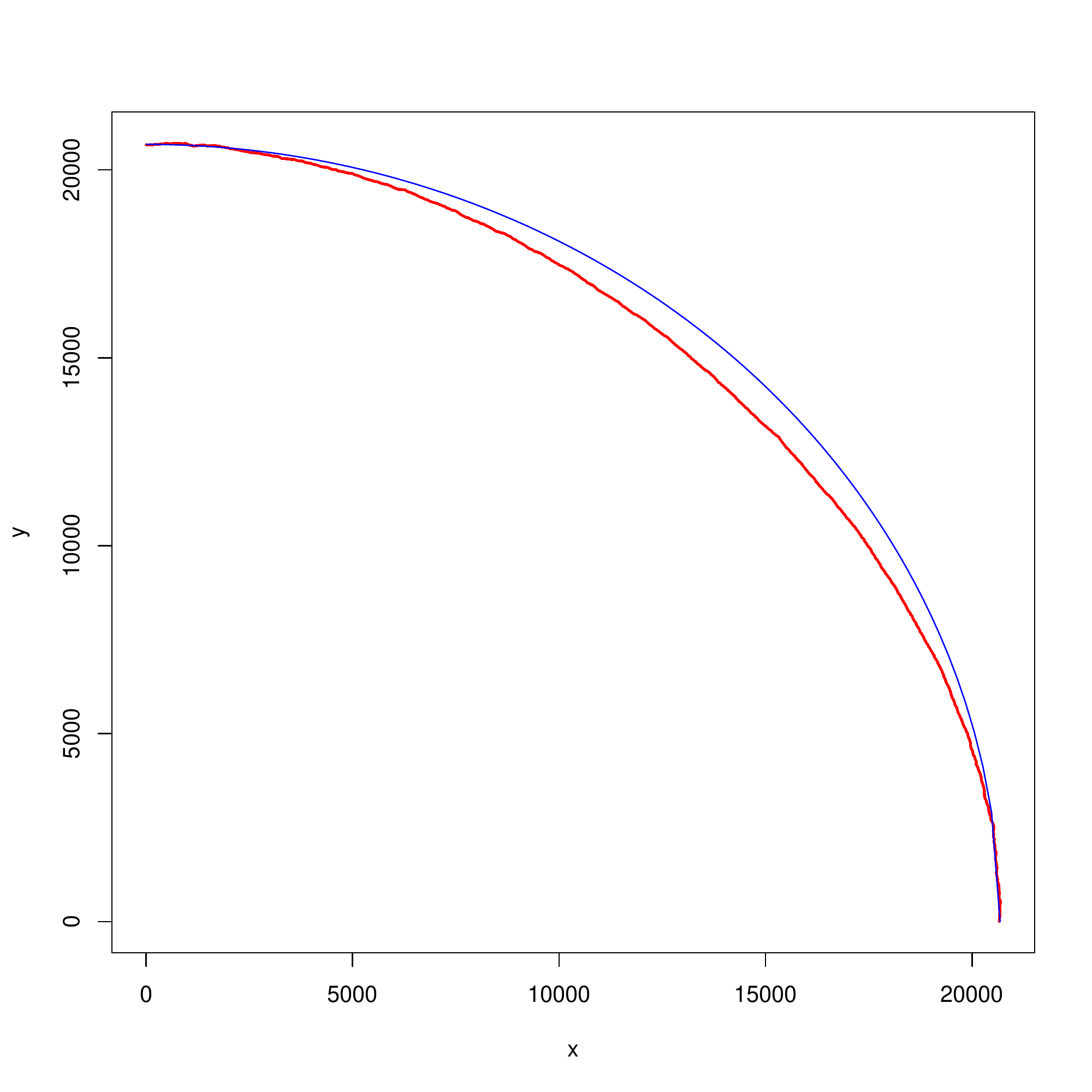}}}\par
\centering \mbox{\subfigure[$k=3$]{\includegraphics[height=7.5cm]{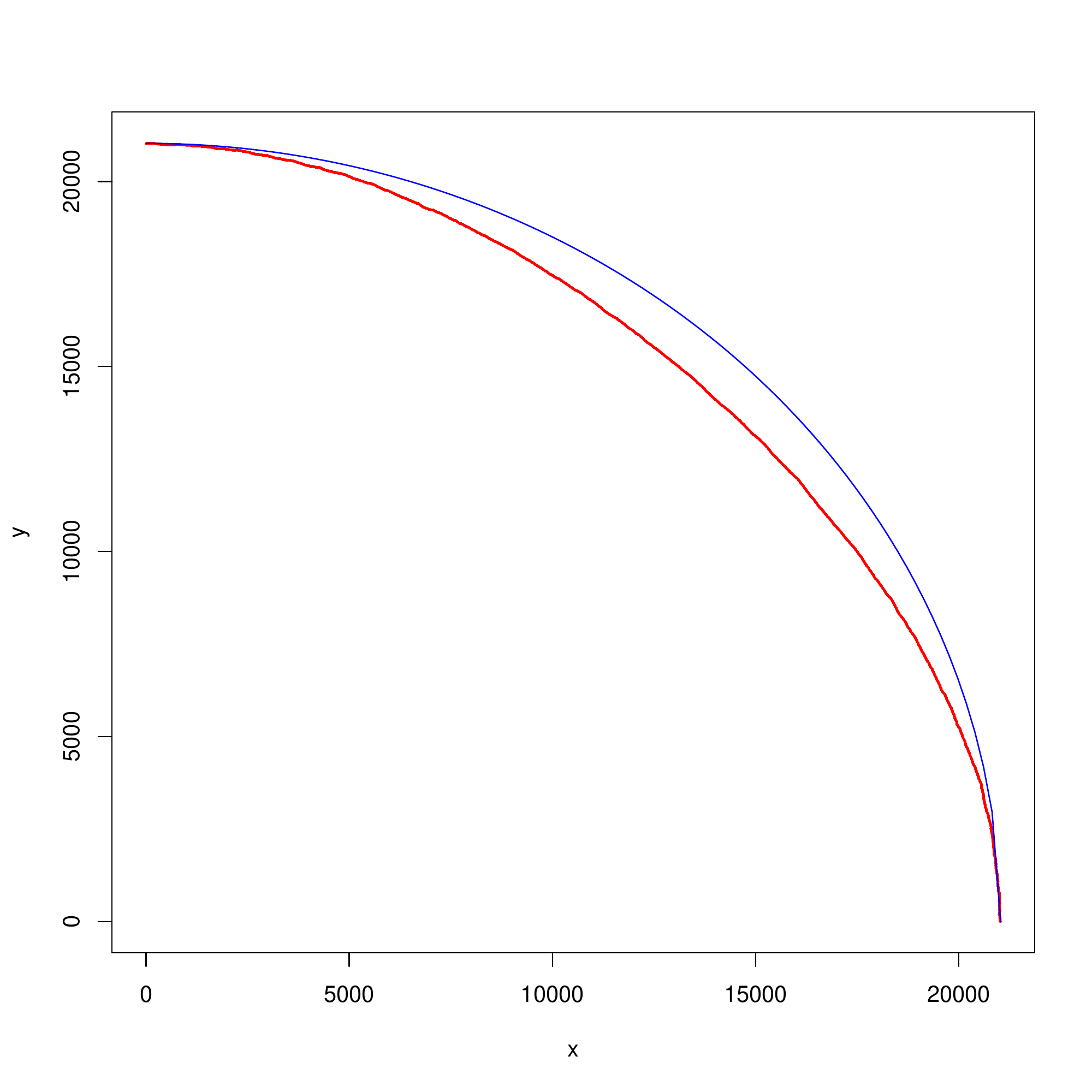}}}\quad
\mbox{\subfigure[$k=4$]{\includegraphics[height=7.5cm]{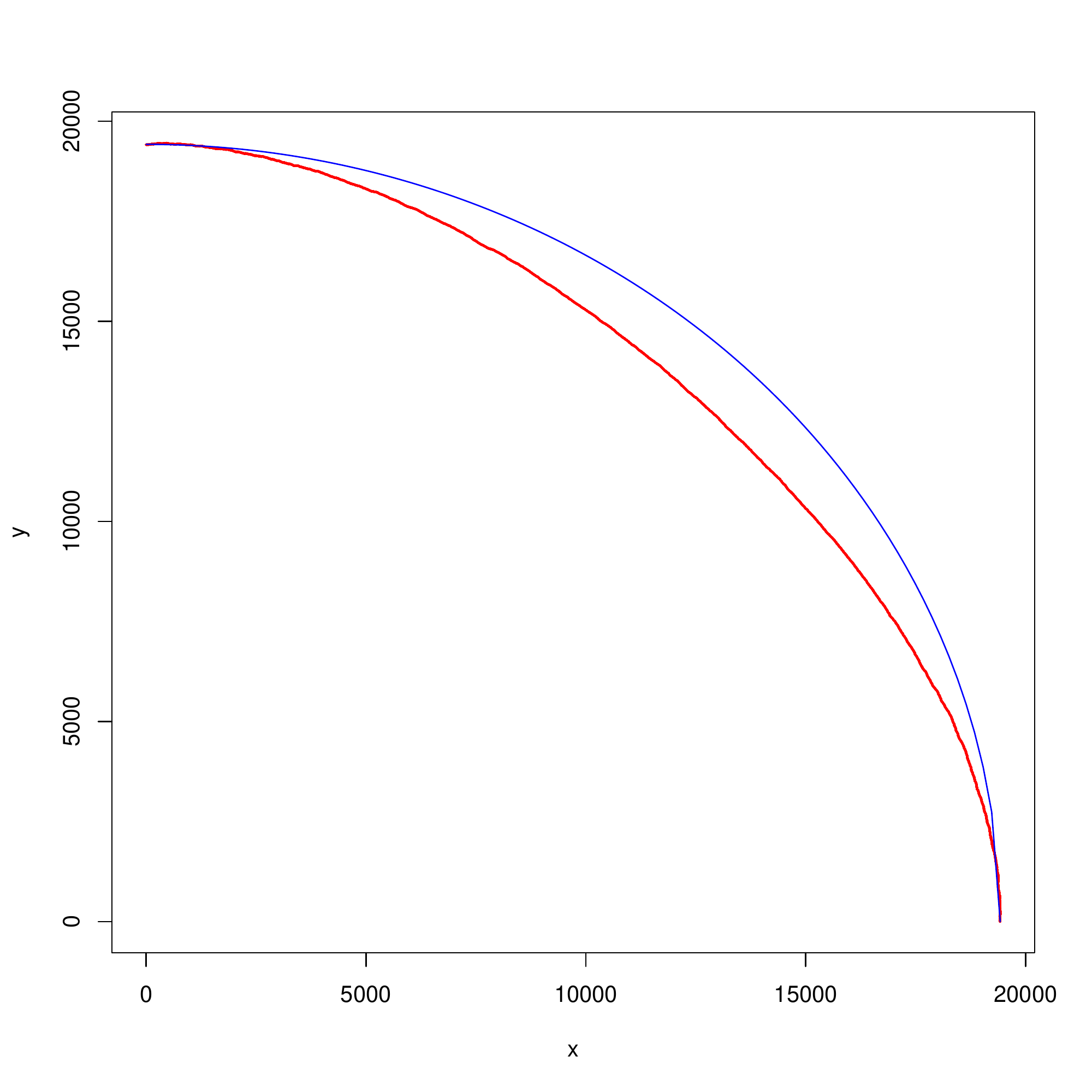}}}\par\caption{Boundary of the infected set for $\Gamma(k,k)$ distributions.}\label{fig:gamma_shape}
\end{figure}

\begin{figure}
\centering \mbox{\subfigure[$\gamma=1$]{\includegraphics[height=7.5cm]{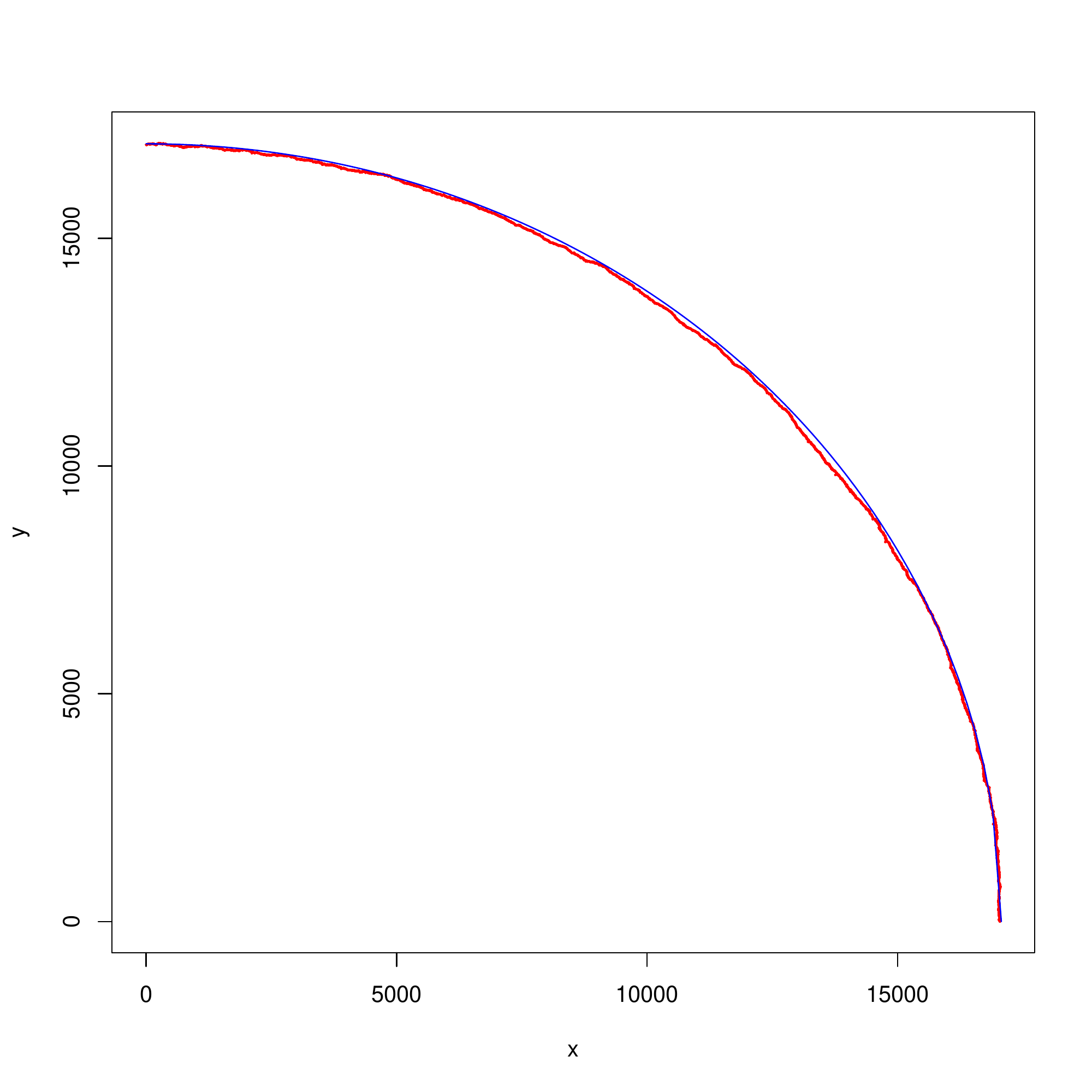}}}\quad
\mbox{\subfigure[$\gamma=2$]{\includegraphics[height=7.5cm]{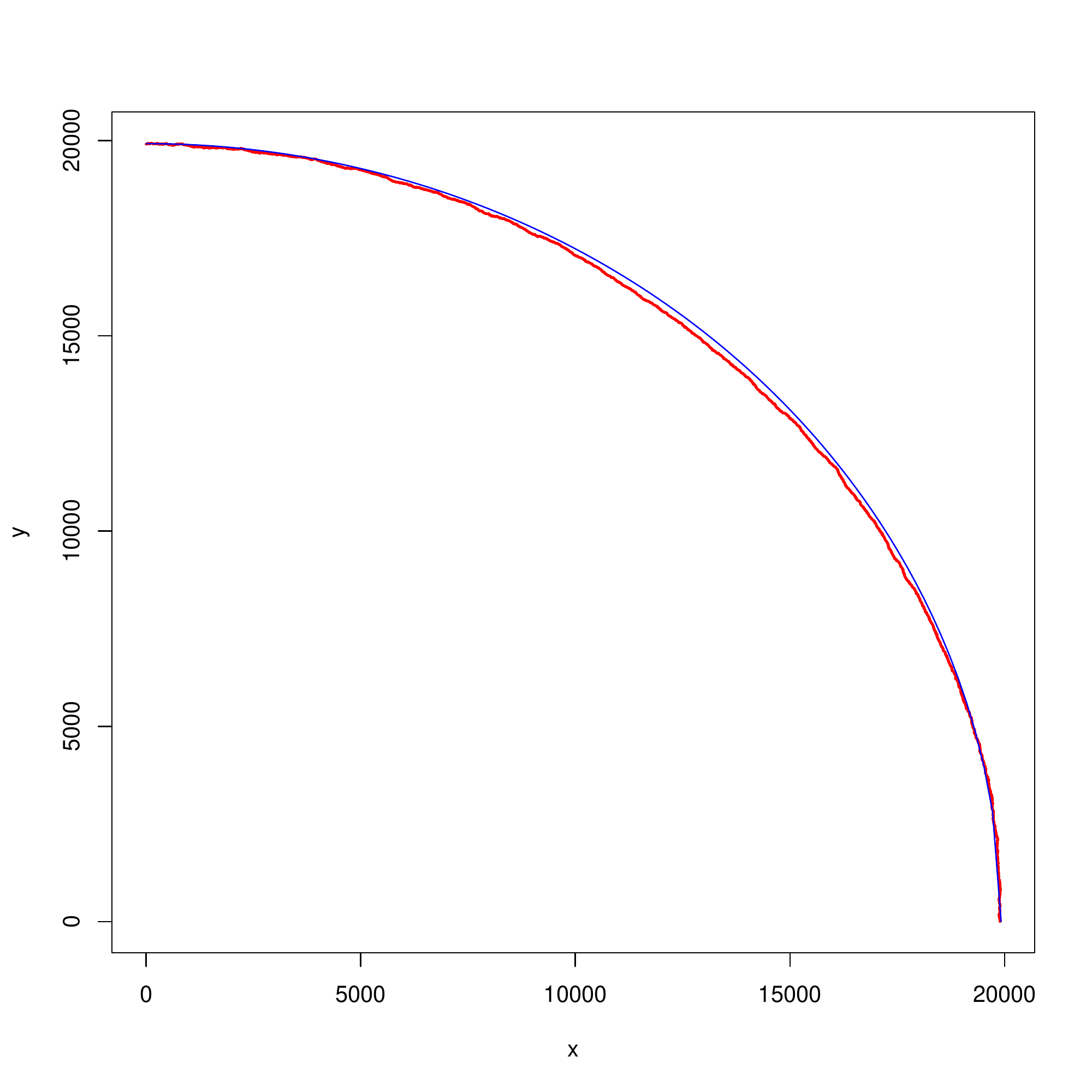}}}\par
\centering \mbox{\subfigure[$\gamma=3$]{\includegraphics[height=7.5cm]{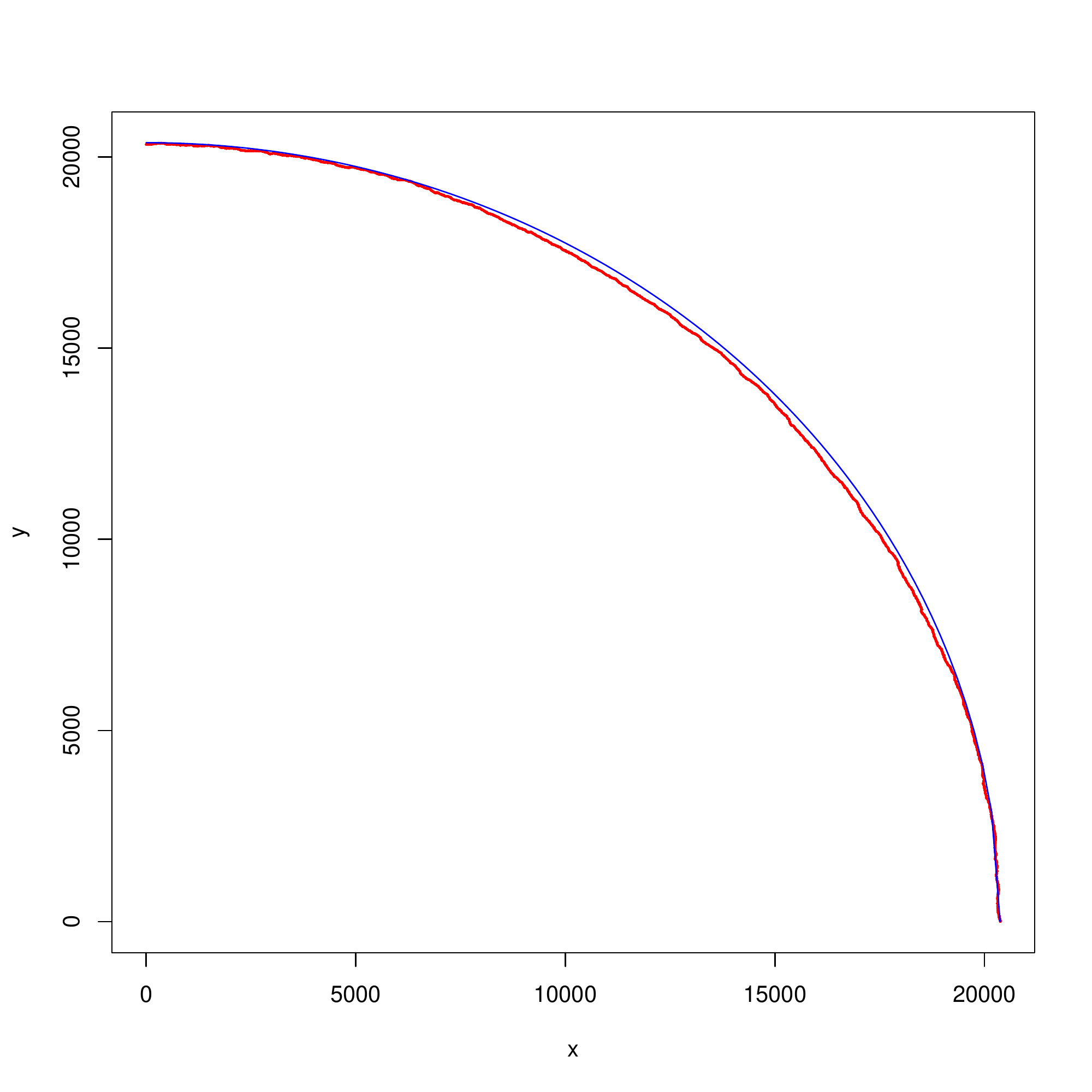}}}\quad
\mbox{\subfigure[$\gamma=4$]{\includegraphics[height=7.5cm]{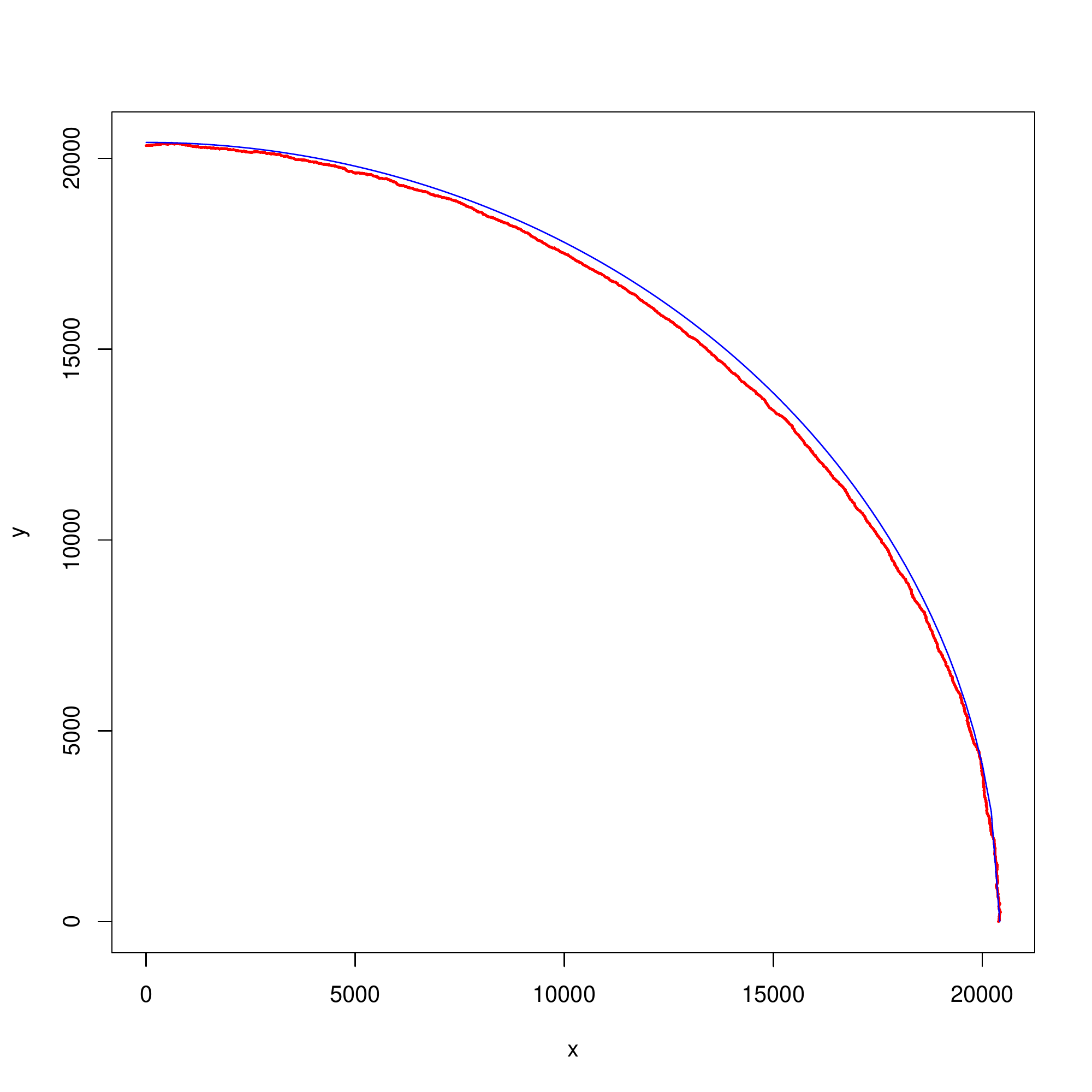}}}\par\caption{Boundary of the infected set for $\widehat{\rm{Fi}}(\gamma)$ distributions.}\label{fig:Fi_shape}
\end{figure}

\begin{figure}
\centering
\begin{minipage}{.48\textwidth}
  \begin{center}
  \includegraphics[width=6cm]{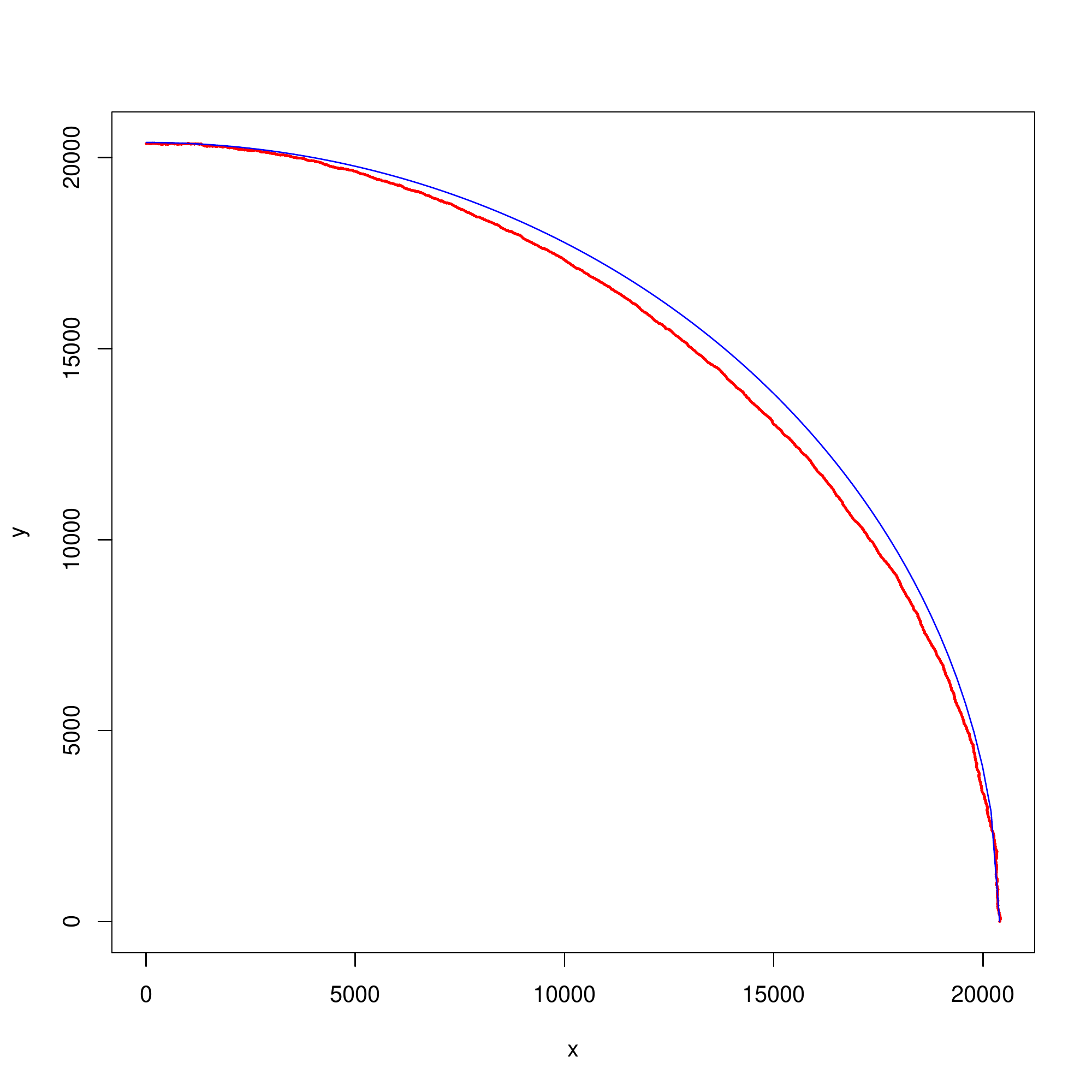}
  \end{center}
  \caption{Boundary of the infected set for U(0,1).}
  \label{fig:U_shape}
\end{minipage}
\hspace{0.2cm}
\begin{minipage}{.48\textwidth}
  \begin{center}
  \includegraphics[width=7.5cm]{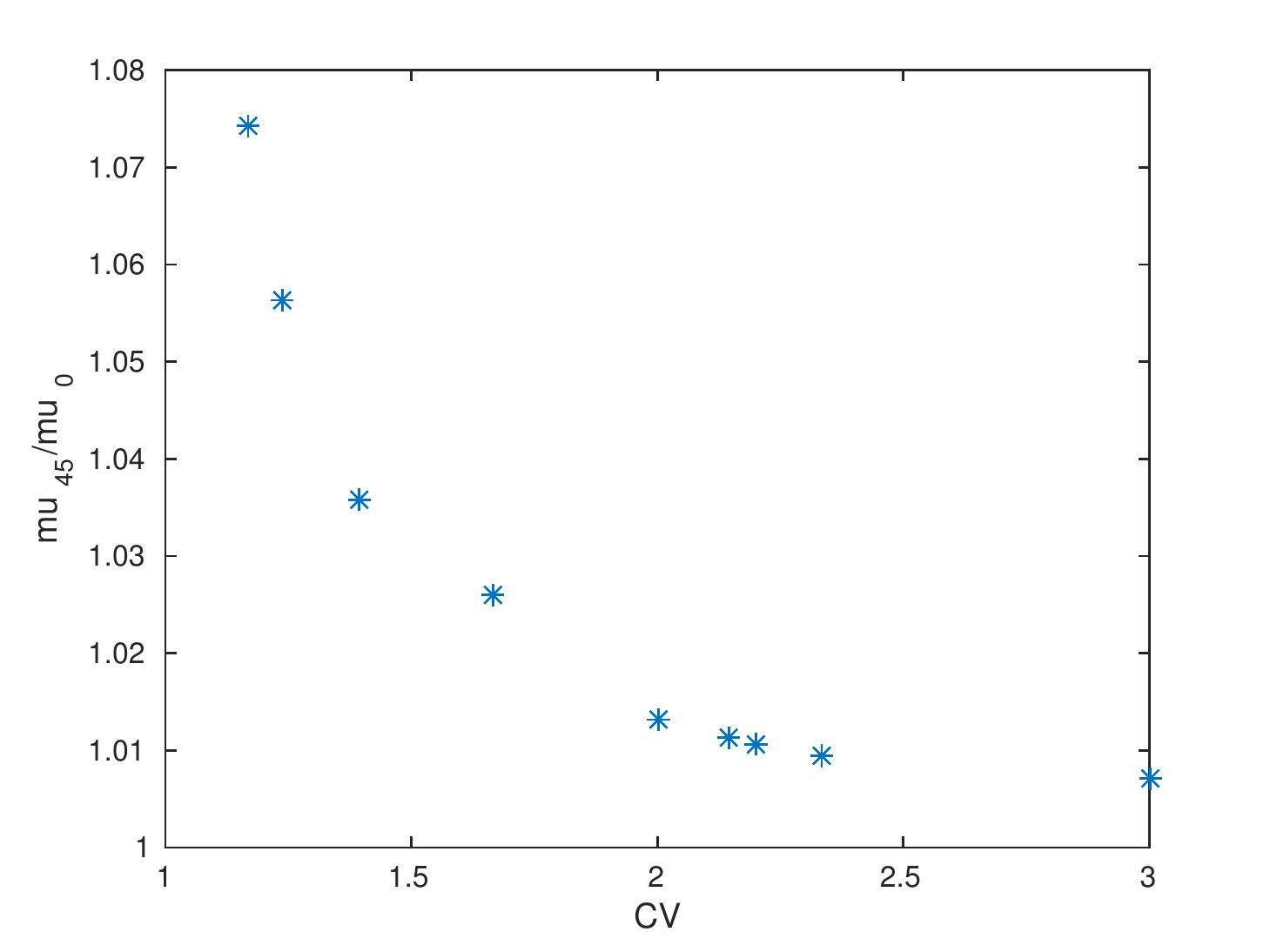}
  \end{center}
  \caption{The estimated ratio $\mu_{45}/\mu_0$ plotted against $\E_4[\tau^2]/\E_4^2[\tau]$.}
  \label{fig:kvot_CV}
\end{minipage}
\end{figure}

We now turn to the asymptotic shape. The last paragraph in the previous section indicates that, for all simulated distributions, the asymptotic shape is contained in a circle with radius given by the speed $\mu^{-1}$ of the growth along the axes. Figure \ref{fig:gamma_shape} shows the boundary of the infected sets (in single realizations) for $\Gamma(k,k)$-distributions with $k=1,2,3,4$. We have first checked roughly how long the simulations have to run to reach the point $(20.000,0)$ and the pictures are then based on realizations up to such fixed times (which hence may be different for different pictures). A circle with radius given by the estimated asymptotic speed of the growth along the axes is included for comparison. Figure \ref{fig:Fi_shape} shows corresponding pictures for $\widehat{\rm{Fi}}(\gamma)$ distributions with $\gamma=1,2,3,4$, and Figure \ref{fig:U_shape} for the uniform distribution U(0,1). From the pictures we see that the shape seems to be closer to a circle for distributions with larger variability. Indeed, for the Gamma distributions, the variance increases as $k$ decreases, the $\widehat{\rm{Fi}}$ distributions are more heavy-tailed for smaller values of $\gamma$. Comparing the different distribution types, we also see that the shapes for the heavy-tailed Fisher type distributions are the most circle-like. 

One possible quantification of the deviation from a circle is the ratio $\mu_{45}/\mu_0$ -- if $\mu_{\theta}$ is indeed increasing as $\theta\in[0,45]$ increases, the difference between the shape and a circle is maximal for $\theta=45$. In Figure \ref{fig:kvot_CV}, estimated values of these ratios are plotted against the squared coefficients of variation CV$_4^2=\E_4[\tau^2]/\E_4^2[\tau]$ for the distributions in Table \ref{tab_summary} (recall that $\E_4[\tau^2]<\infty$ is a necessary condition for convergence to an asymptotic shape). Here, CV$_4$ measures the standard deviation of $\min\{\tau_1,\ldots,\tau_4\}$ as a fraction of the mean. We see that there indeed seems to be a relation in that a larger value of $\rm{CV}^2$ tends to imply a smaller value of $\mu_{45}/\mu_0$.

\begin{figure}
\centering \mbox{\subfigure[Exp(1)]{\includegraphics[height=5cm]{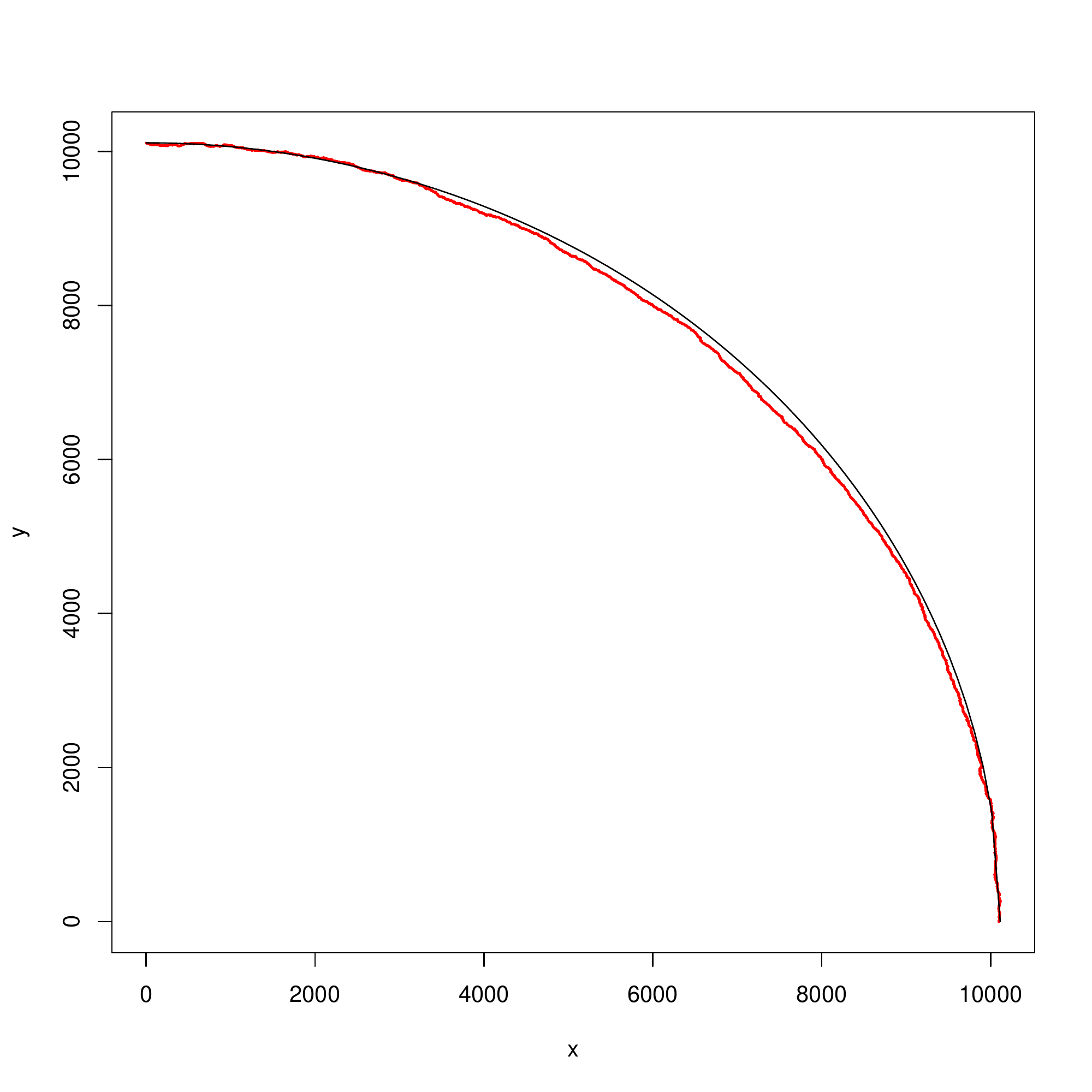}}}\quad
\mbox{\subfigure[0.5+Exp(1)]{\includegraphics[height=5cm]{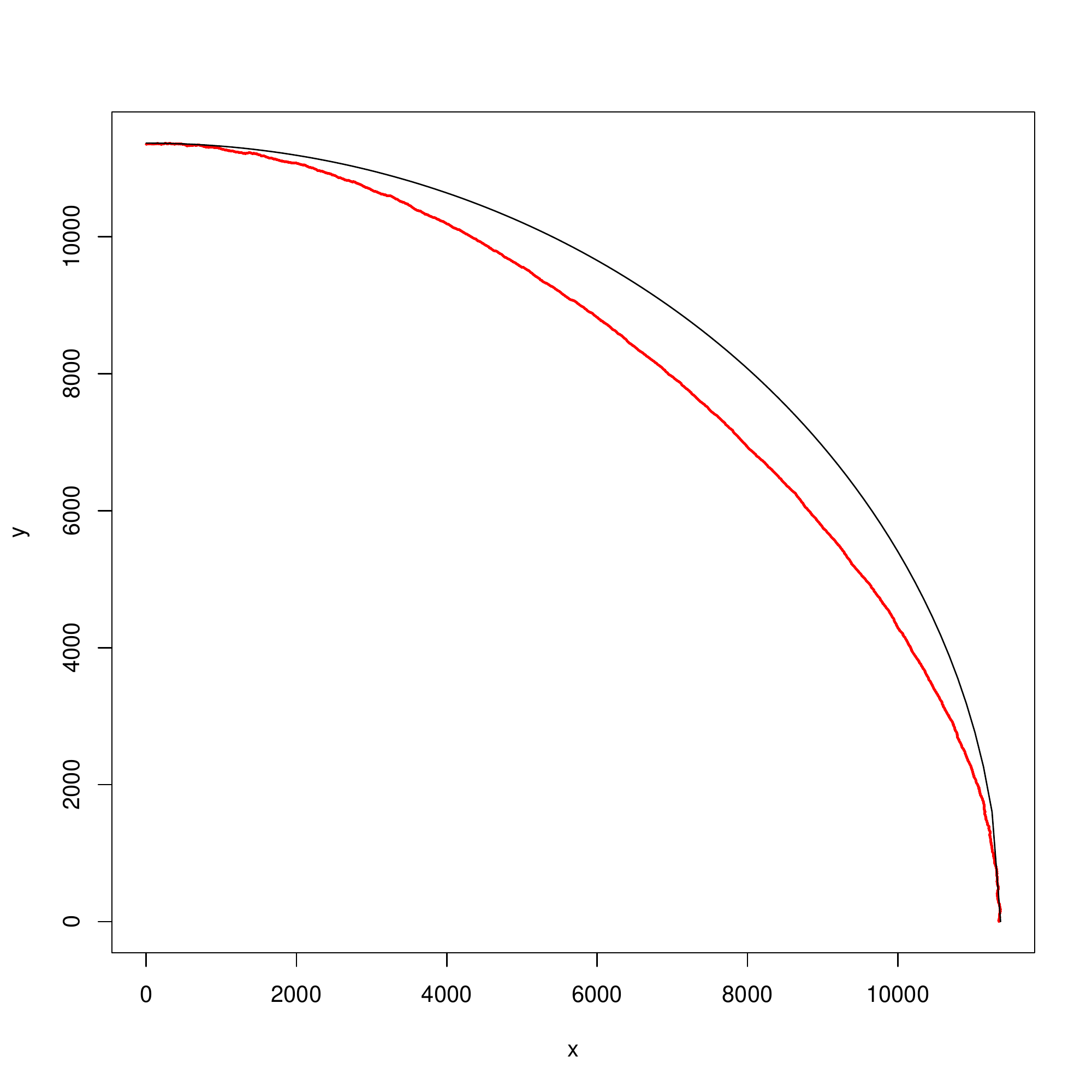}}}\quad\mbox{\subfigure[4+Exp(1)]{\includegraphics[height=5cm]{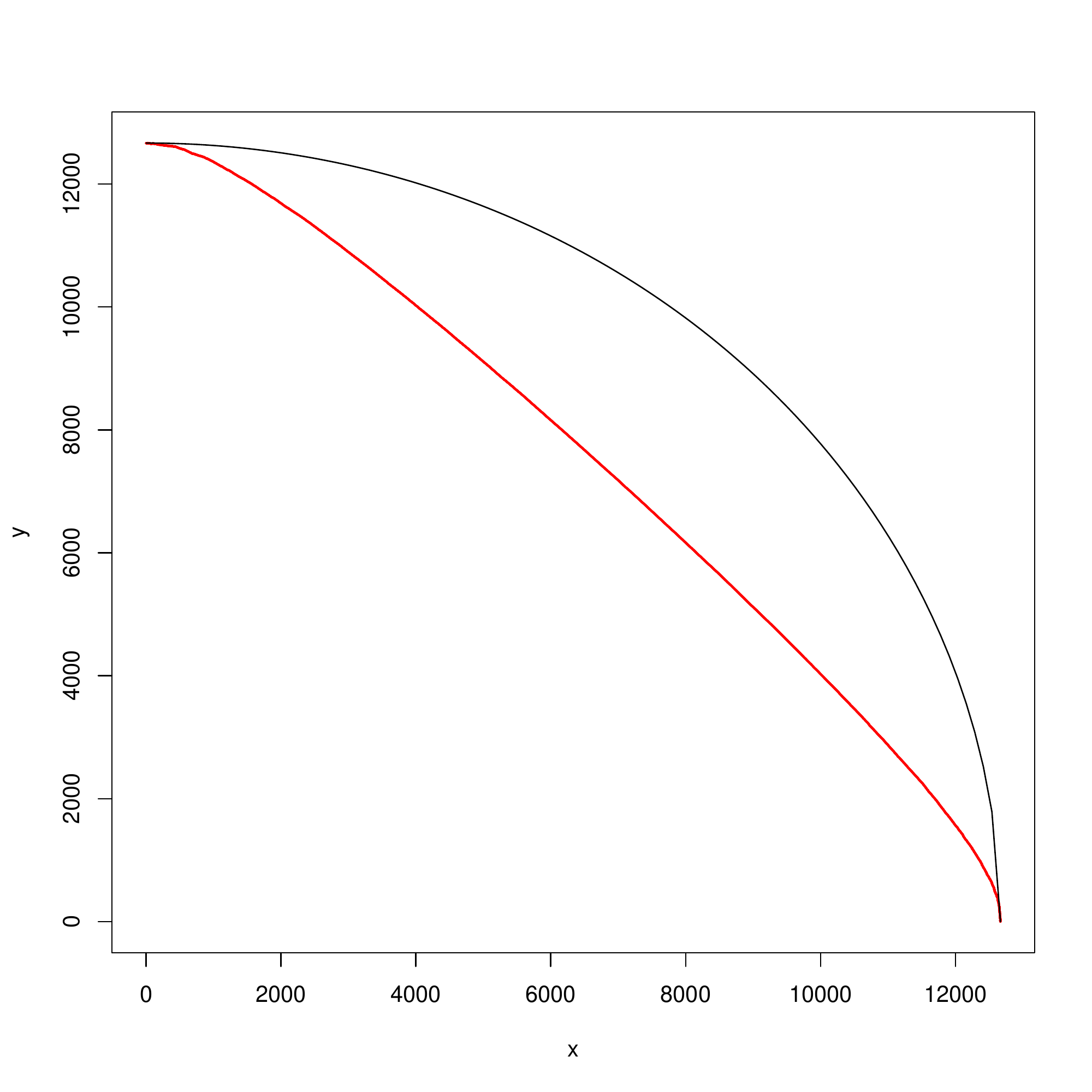}}}\par\caption{Boundary of the infected set for translated exponential distribution.}\label{fig:exptrans_shape}
\end{figure}

\begin{figure}
\centering \mbox{\subfigure[$\widehat{\rm{Fi}}(2)$]{\includegraphics[height=5cm]{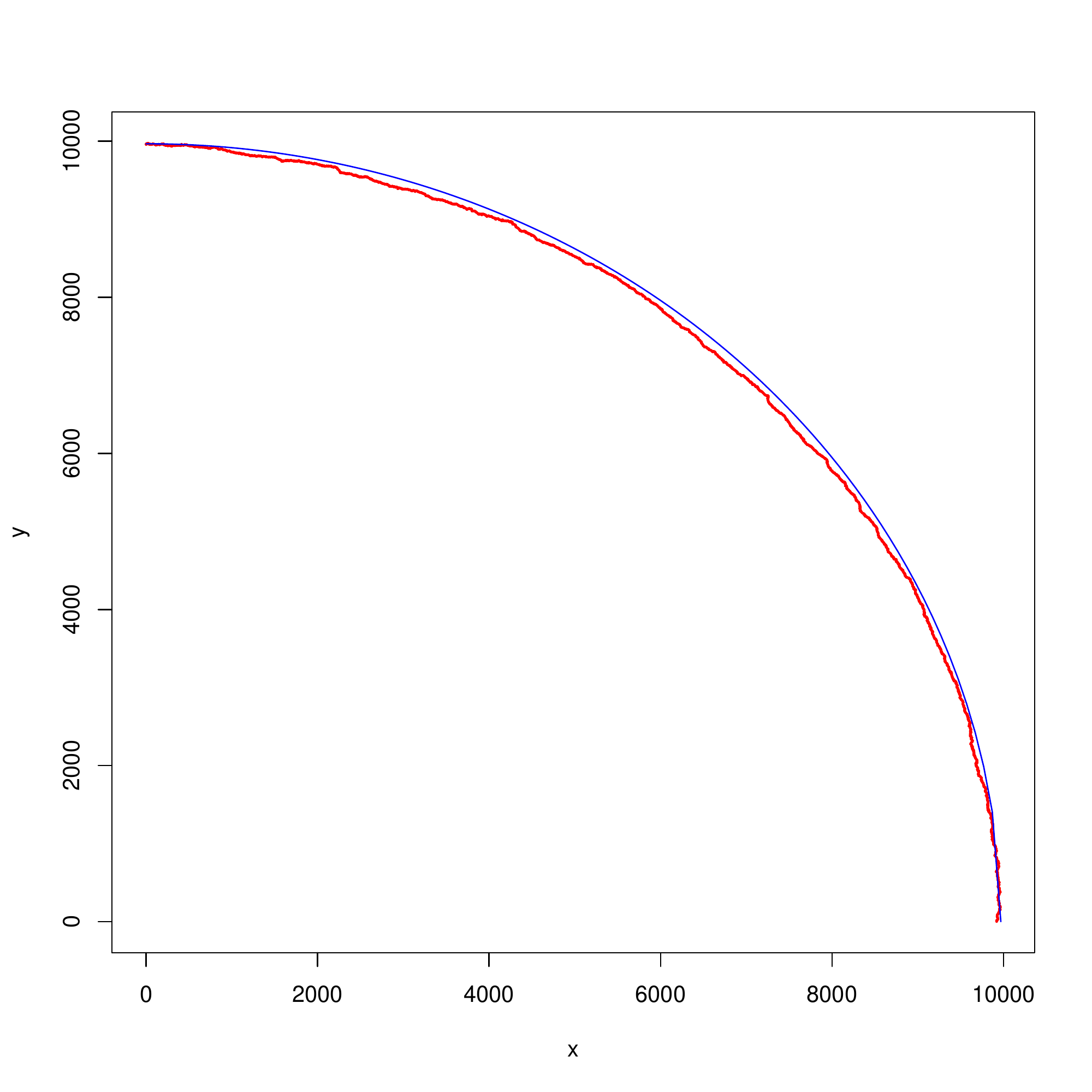}}}\quad
\mbox{\subfigure[$0.5+\widehat{\rm{Fi}}(2)$]{\includegraphics[height=5cm]{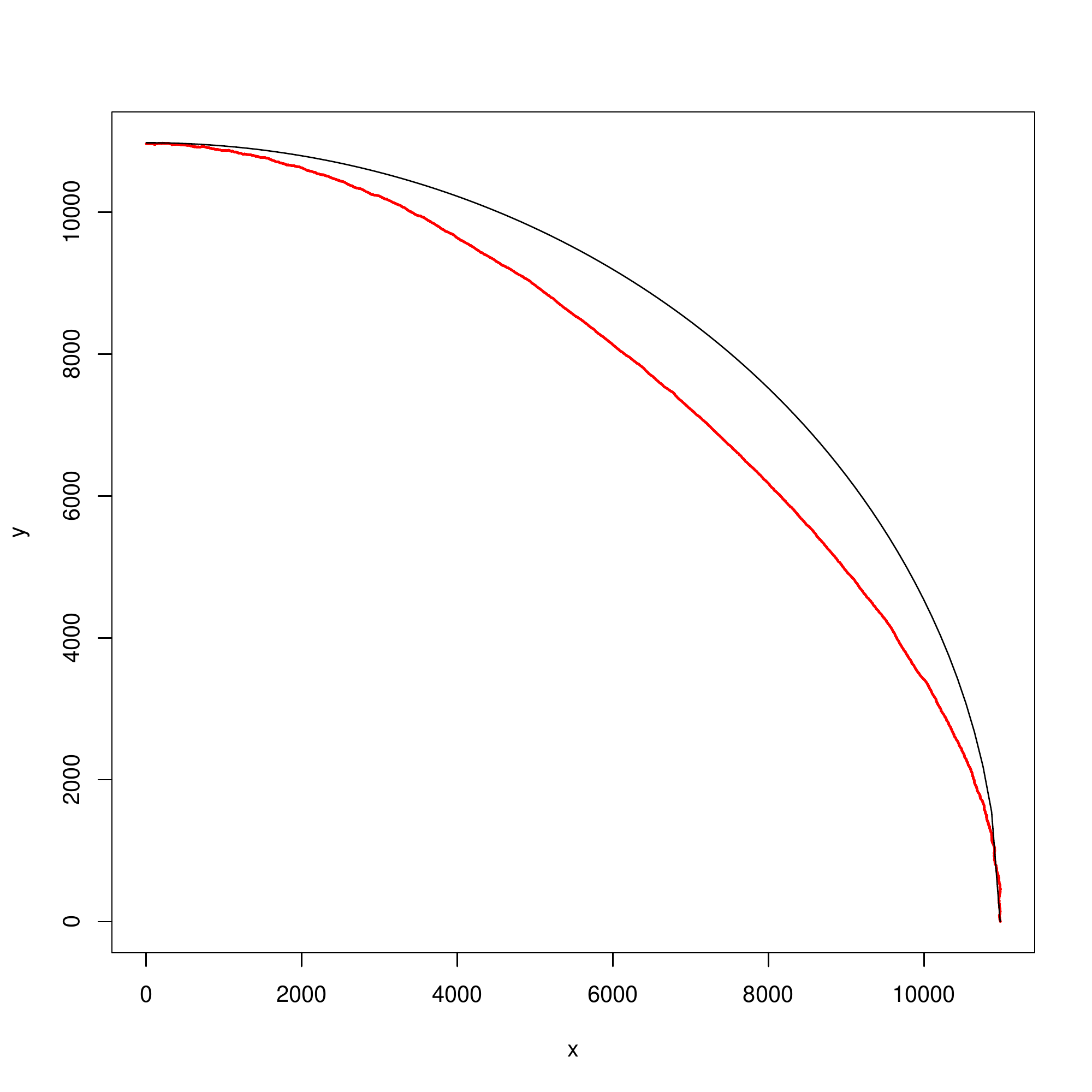}}}\quad\mbox{\subfigure[$4+\widehat{\rm{Fi}}(2)$]{\includegraphics[height=5cm]{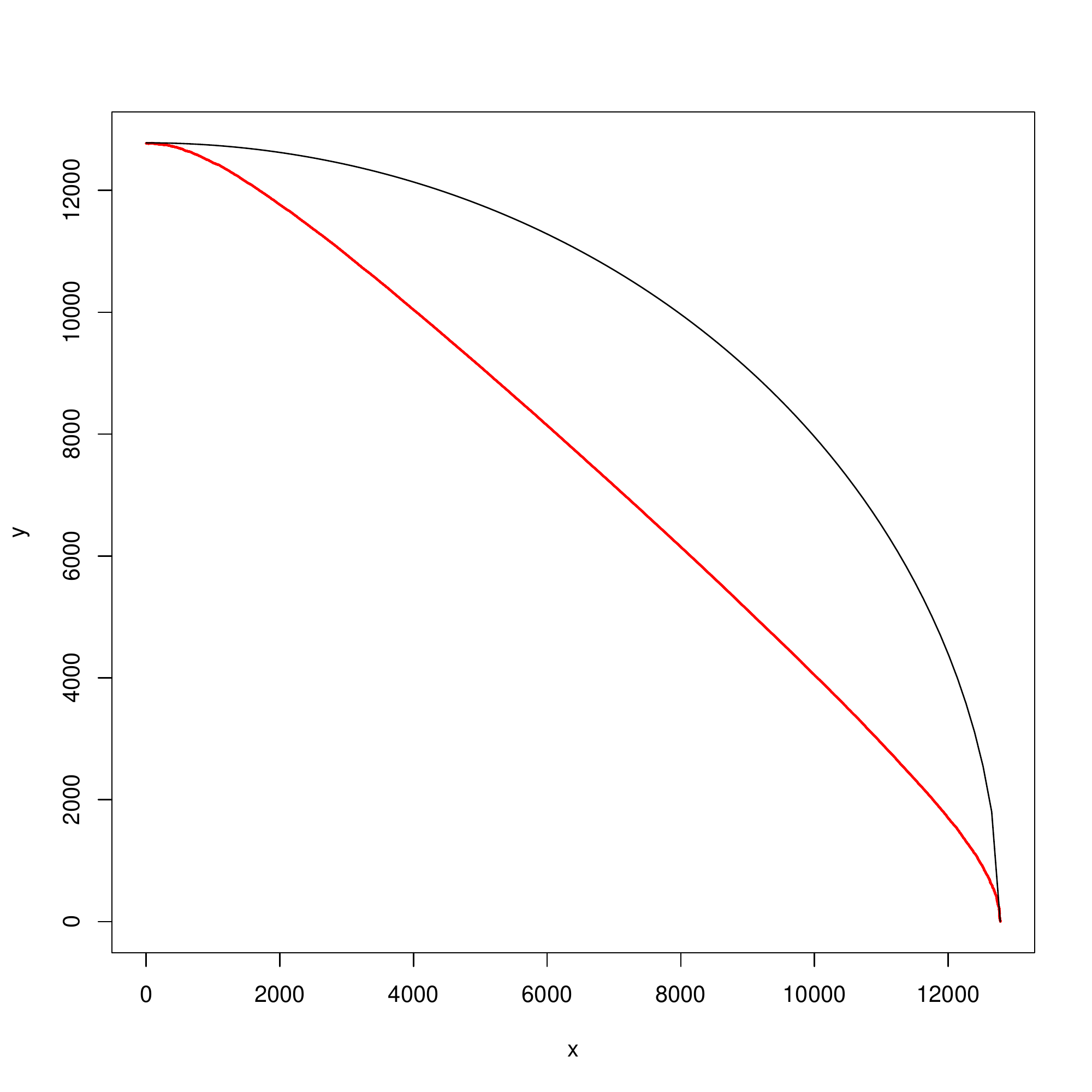}}}\par\caption{Boundary of the infected set for translated $\widehat{\rm{Fi}}(2)$ distribution.}\label{fig:Fitrans_shape}
\end{figure}

\begin{figure}
\centering \mbox{\subfigure[$\Gamma(2,2)$]{\includegraphics[height=5cm]{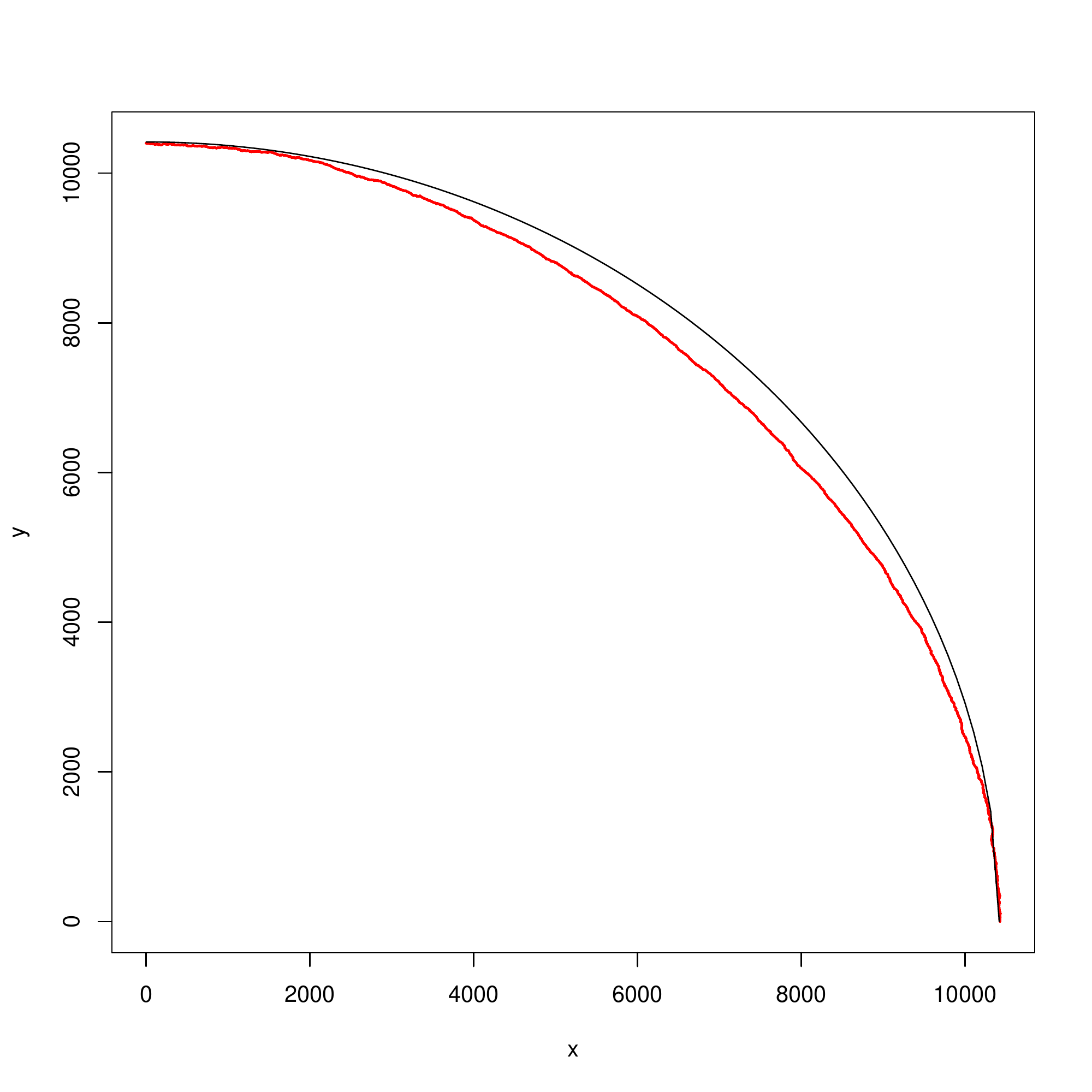}}}\quad
\mbox{\subfigure[$\Gamma(10,10)$]{\includegraphics[height=5cm]{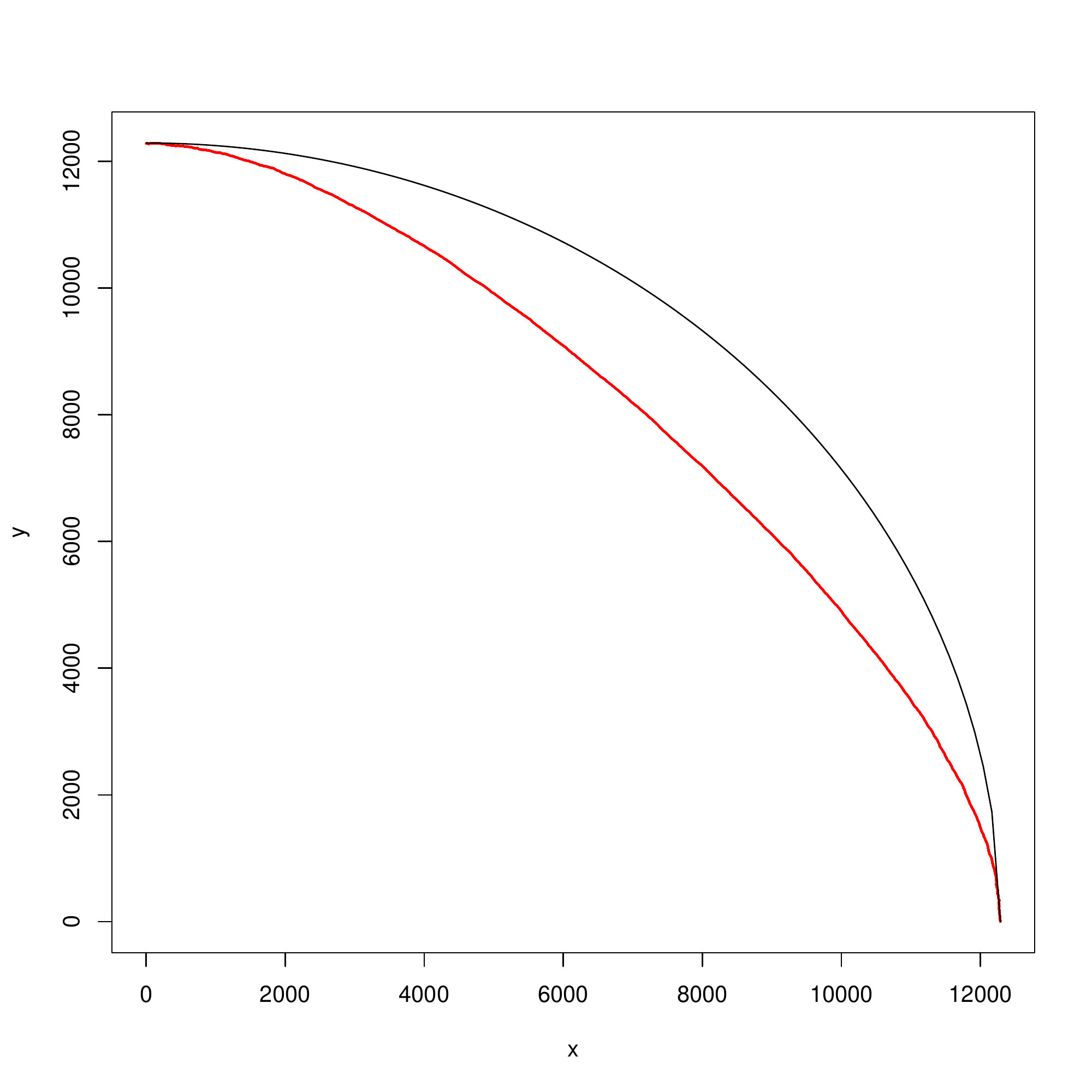}}}\quad\mbox{\subfigure[$\Gamma(20,20)$]{\includegraphics[height=5cm]{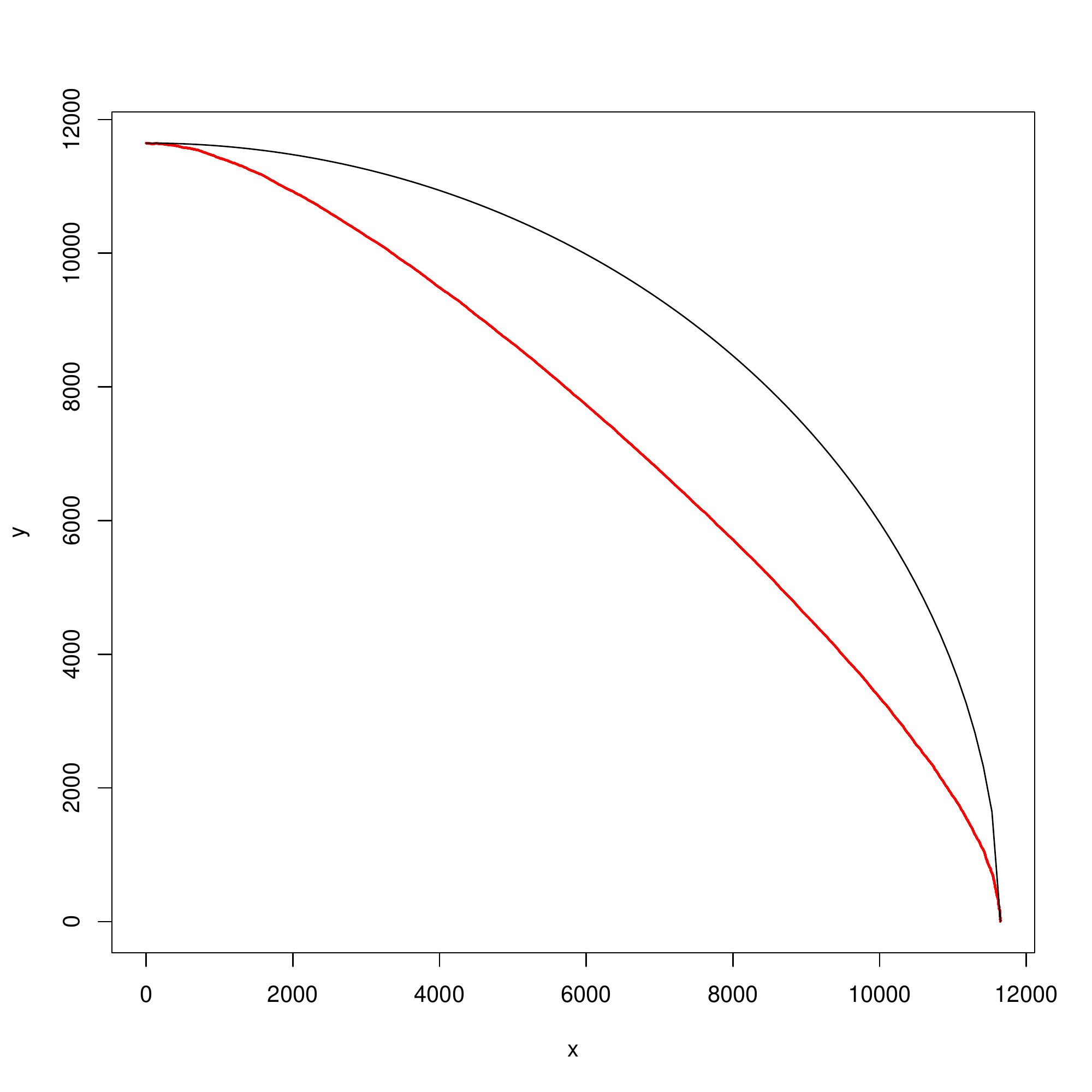}}}\par\caption{Boundary of the infected set for Gamma distributions with decreasing variance.}\label{fig:gammatrans_shape}
\end{figure}

Figures \ref{fig:exptrans_shape} and \ref{fig:Fitrans_shape} show pictures of the shape for the translated exponential and $\widehat{\rm{Fi}}(2)$ distributions, with pictures for the non-translated distributions included for comparison (also with $x_{max}=10.000$). As described above, translating the passage time distribution brings the model closer to a situation with constant passage times and indeed we see that the shape comes closer to a diamond -- the asymptotic shape for constant passage times -- as the translation increases. The pictures are similar for other translated distributions. Figure \ref{fig:gammatrans_shape} shows that, when the parameter increases in the Gamma distribution, so that the variance decreases, then the shape also comes closer to a diamond.

\begin{figure}
\centering \mbox{\subfigure[$\widehat{\rm{Fi}}(0.5)$]{\includegraphics[height=7.5cm]{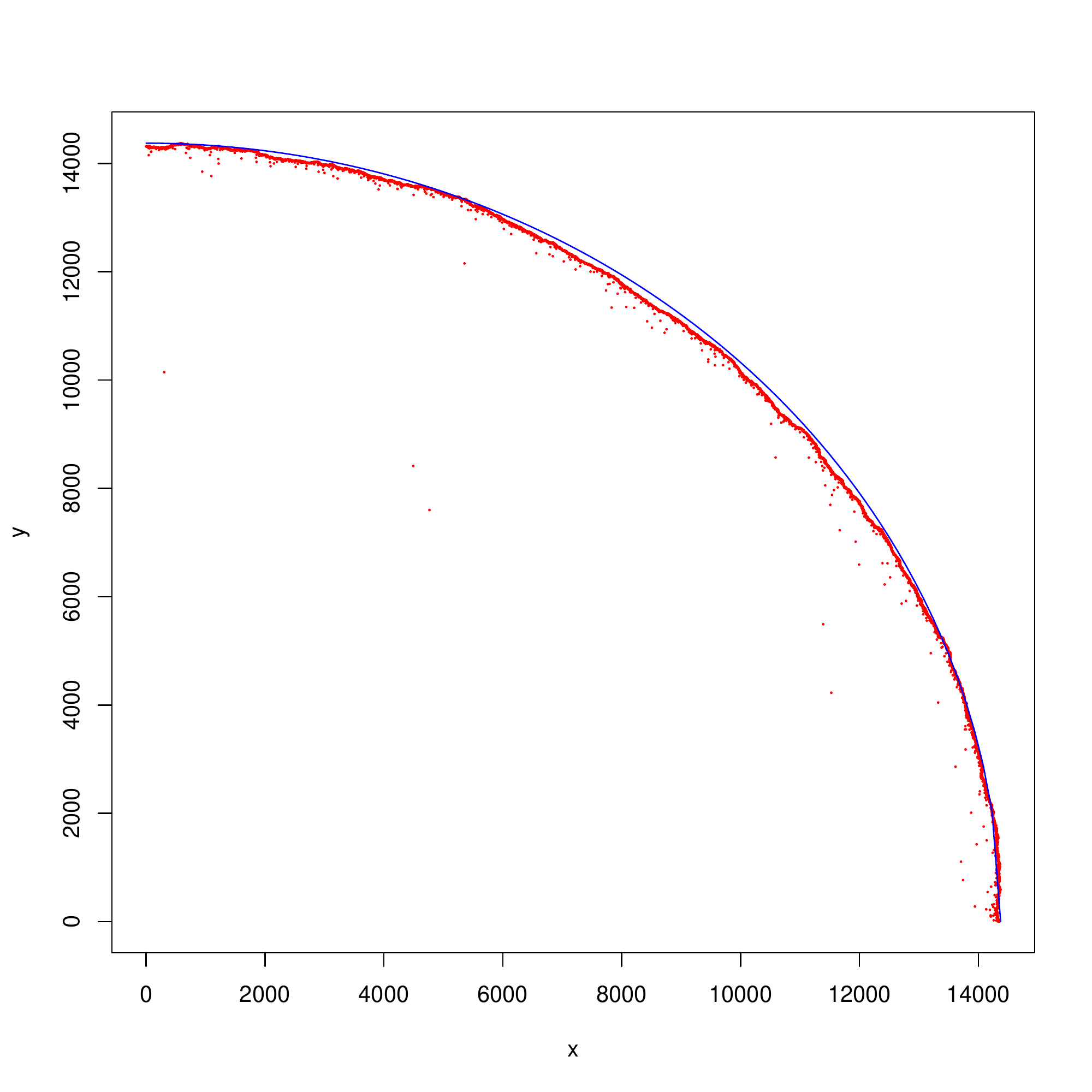}}}\quad
\mbox{\subfigure[$\widehat{\rm{Fi}}(0.3)$]{\includegraphics[height=7.5cm]{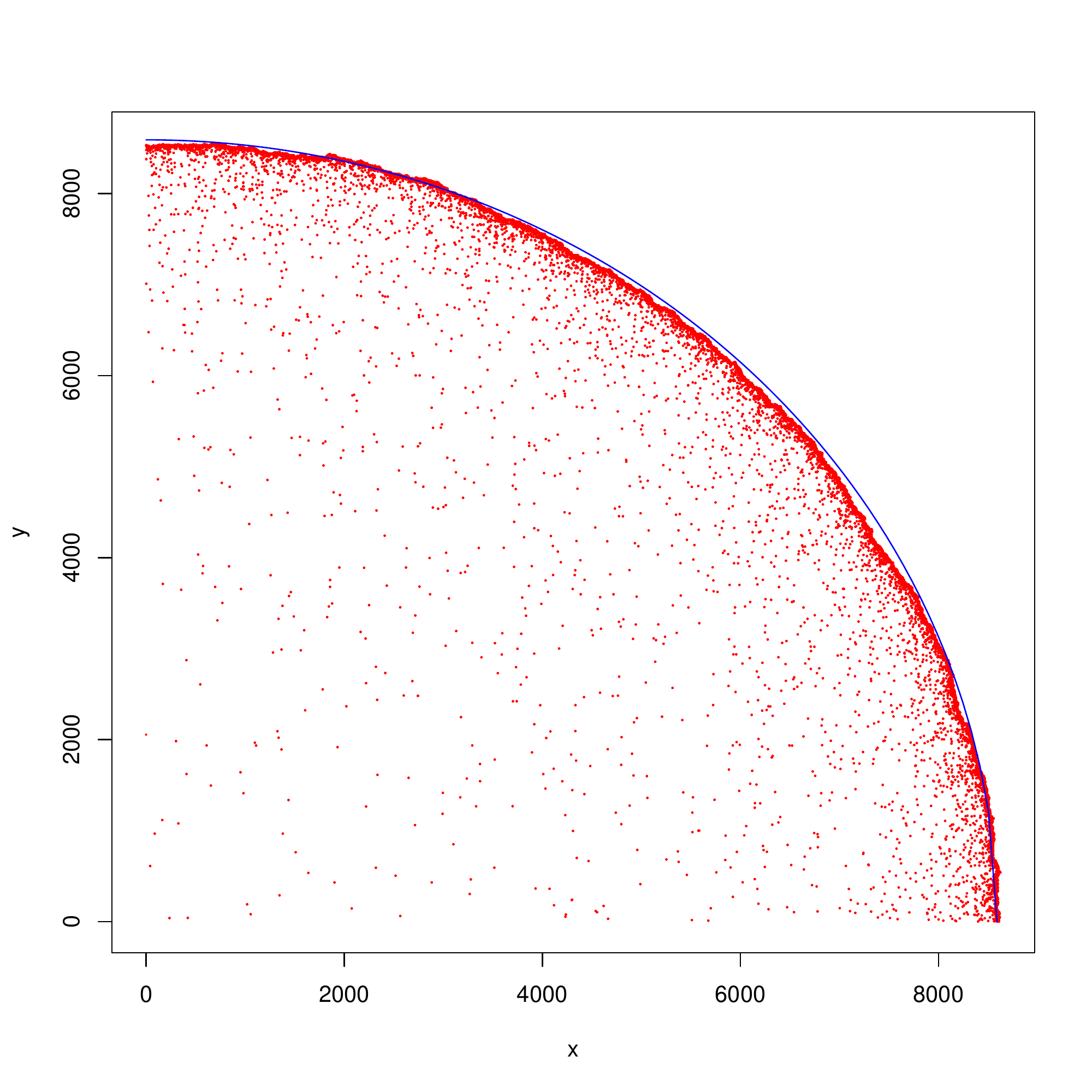}}}\par\caption{Boundary of the infected set for heavy tailed Fisher distributions.}\label{fig:Fi_heavy_shape}
\end{figure}

When $\E_4[\tau^2]=\infty$, the conclusions of the shape theorem are not true. As pointed out after Theorem \ref{th:shape}, the infected region will then have uninfected points in its interior that are visible on a linear scale. This is illustrated in Figure \ref{fig:Fi_heavy_shape}, with pictures of the boundary of the infected set for $\widehat{\rm{Fi}}(0.5)$ and $\widehat{\rm{Fi}}(0.3)$. Note that we still have $\E_4[\tau]<\infty$, so that the time constants do exist.

\begin{figure}[t]
\centering \mbox{\subfigure[Fluctuations of hitting points on lines $L_n$.]{\includegraphics[height=7.5cm]{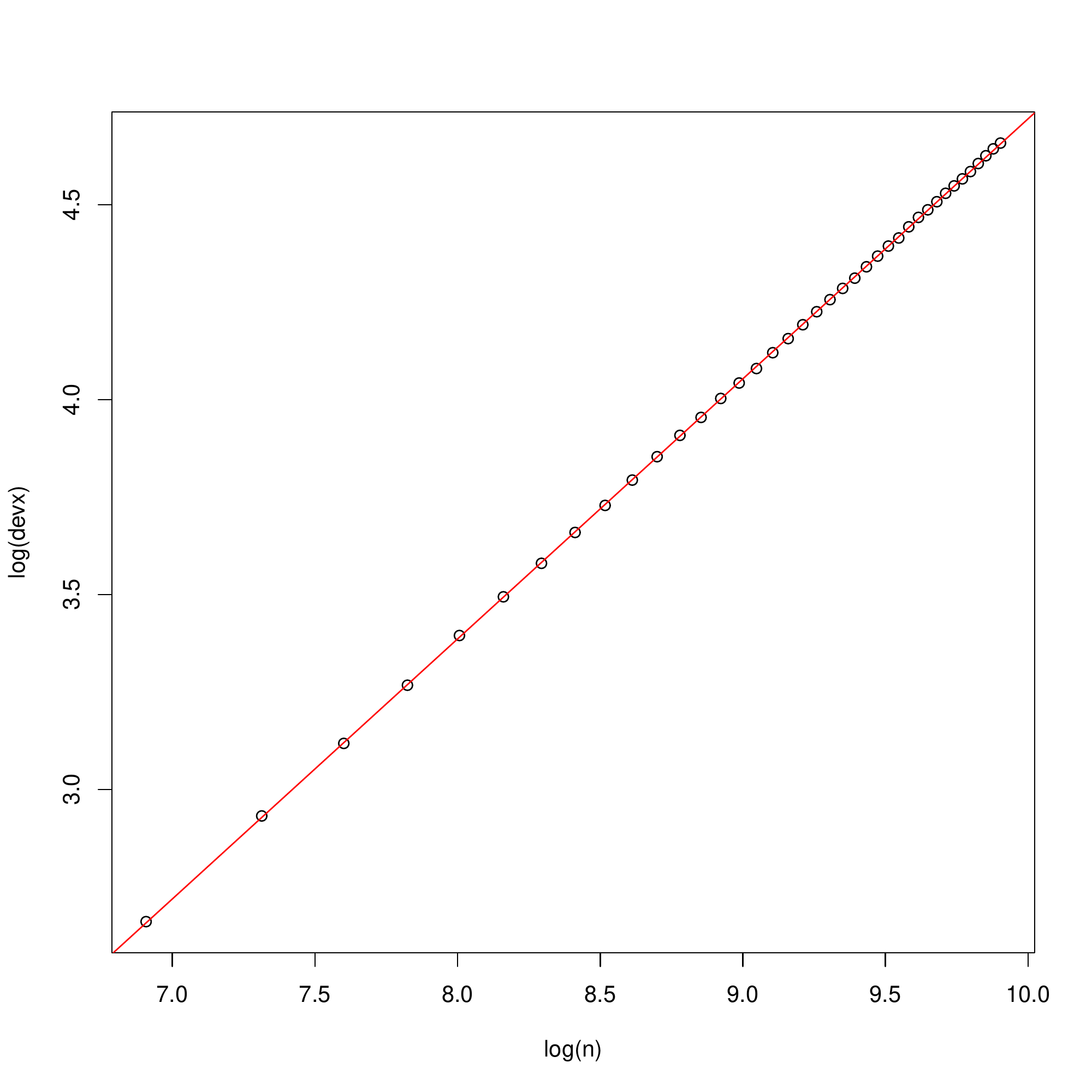}}}\quad
\mbox{\subfigure[Standard deviations of hitting times.]{\includegraphics[height=7.5cm]{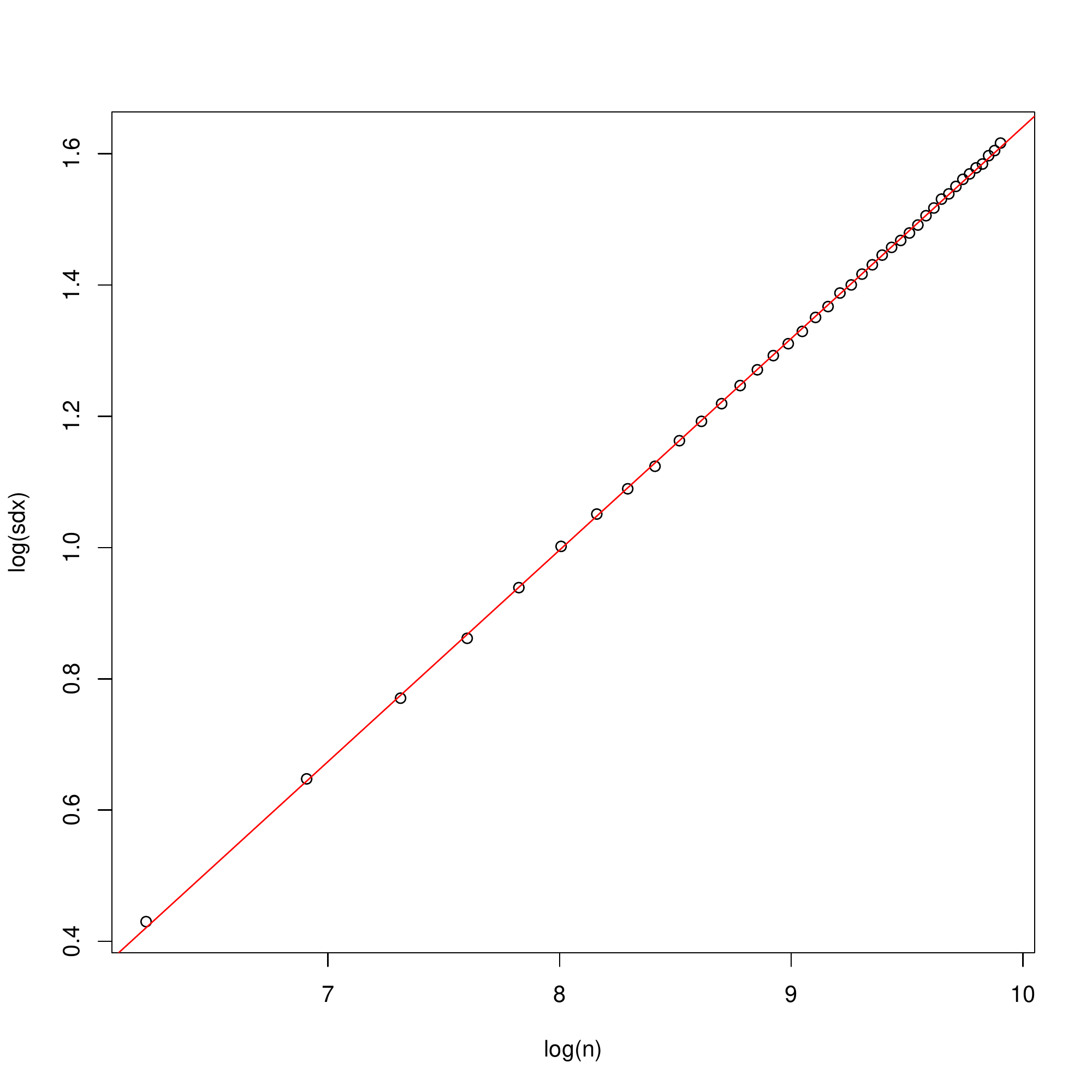}}}\par\caption{Regressions for the exponents $\chi$ and $\xi$ for exponential passage times.}\label{fig:Bscaling}
\end{figure}

Recall from Section \ref{sec:AS} that the fluctuations of the interface between the infected set $\bar{B}(t)$ and its complement is believed to be of the order $t^{\chi}$ with $\chi=1/3$, so that the standard deviation of the hitting time $T(\mathbf{0}, \mathbf{n})$ is of the order $n^{1/3}$, and the fluctuations of the location where $\bar{B}(t)$ hits a line at distance $n$ from the origin is conjectured to be of the order $n^{\xi}$ with $\xi=2/3$. We close this section by referring to Figure \ref{fig:Bscaling} for some support of this. Let $H_n$ denote the distance to the $x$-axis from the hitting point on the vertical line $L_n$ at distance $n$ from the origin. Fitting a regression model with $\E[H_n]=cn^{\xi}$ based on sample points at distance 500, gives $\hat{\xi}=0.6670786$ with standard deviation 0.0005222. Due to the observation points being positively correlated, the standard deviation is likely to be too small, but the predicted value 2/3 is still contained in the 95$\%$ confidence interval. Figure \ref{fig:Bscaling}(a) shows a plot of $\log(H_n)$ against $\log n$ and we see that the observations indeed fit very well with the regression line. Similarly we have fitted an analogous regression model for the standard deviation of $T(\mathbf{0},L_n)$, which gives $\hat{\chi}=0.3223056$ with standard deviation 0.0006056. Although the estimate is fairly close to the predicted 1/3, it is somewhat smaller and the 95$\%$ confidence interval does not contain 1/3. Recall however that the standard deviation is likely to be too small so that the confidence interval shold in fact be wider. A contributing factor may also be that we have recorded hitting times to lines rather than to points and such hitting times of course have smaller variance. Figure \ref{fig:Bscaling}(b) shows the observations along with the regression line.

\section{Summary}

We have studied the speed and shape of first passage percolation on $\mathbb{Z}^2$ by aid of computer simulations. Focusing on continuous passage time distributions, the simulated distributions include the Exponential distribution, Gamma distributions, uniform distributions and power-law distributions. We refer to Table \ref{tab_summary} for numerical estimates of time constants and to Section \ref{sec:pictures} for simulation pictures. 

For the time constant, we have found that the most important quantity for determining its value is $\E_4[\tau]$. Indeed, Figure \ref{fig:tc_E4} reveals a close to linear relation for the simulated distributions. For translated versions of the distributions, the relation is less clear along the axis, but becomes more linear closer to the diagonal. For the asymptotic shape, we have found that, for all simulated distributions, the shape is contained in a circle with radius $\mu_0$. Constant passage times give a diamond as asymptotic shape and, when the variance decreases for a given distribution type, the shape comes closer to the diamond. When the variance increases, on the other hand, the shape approaches the circle. In Figure \ref{fig:kvot_CV}, the deviations from the circles, quantified by the ratio $\mu_{45}/\mu_0$, are plotted against the variational coefficient $\E_4[\tau^2]/\E_4^2[\tau]$ for all simulated distributions and indeed we see a negative relation. Finally, we have investigated the exponents for the standard deviation of hitting times and for the fluctuations of hitting points on lines  and our findings are in line with the predicted values 1/3 and 2/3, respectively.

We finish with a few ideas for theoretical investigations:

\begin{itemize}
\item[1.] As mentioned in Section \ref{sim_tc}, we believe that $\mu\geq \E_4[\tau]$ is true in great generality. Is it possible to find sufficient conditions on the passage time distribution for this?
\item[2.] Prove that $\mu_\theta$ is increasing in $\theta\in[0,45]$, perhaps for some suitable class of passage time distributions. This would imply that the asymptotic shape is indeed contained in the circle. 
\item[3.] Does the shape come arbitrarily close to a circle when the variational coefficient $\E_4[\tau^2]/\E_4^2[\tau]$ (or some other quantification of the variability of the distribution) becomes large? We believe that the answer may be no. 
\end{itemize}

\pagebreak

\section{Simulation pictures}\label{sec:pictures}

\end{document}